\documentclass[12pt]{article}
\usepackage{amsmath}
\usepackage{amssymb}

\newcommand{\eT}{{\mathbf{T}}}

\def\a{\alpha}
\def\b{\beta}
\def\g{\gamma}
\def\G{\Gamma}

\def\la{\lambda}
\def\Lm{\Lambda}
\def\S{\Sigma}
\def\om{\omega}

\def\t{\tau}
\def\ts{\times}
\def\th{\theta}

\def\wh{\widehat}
\def\wt{\widetilde}
\def\ov{\overline}

\def\BC{{\mathbb C}}

\def\clp{{\mathcal P}}
\def\cla{{\mathcal A}}
\def\clb{{\mathcal B}}

\def\im{{\rm Im\ }}

\newtheorem{Pa}{Paper}[section]
\newtheorem{Tm}[Pa]{{\bf Theorem}}
\newtheorem{La}[Pa]{{\bf Lemma}}

\newtheorem{Pn}[Pa]{{\bf Proposition}}

\newcommand{\ands}{\quad\mbox{and}\quad}
\newcommand{\iy}{{\infty}}

\newcommand{\bpr}{{\noindent\textbf{Proof.}\ }}
\newcommand{\epr}{{\hfill $\Box$}}

\title{Krein systems and canonical systems\\ 
on a finite interval: accelerants \\ with a jump discontinuity at the origin \\ and  continuous potentials}

\author{D. Alpay, I. Gohberg (Z$''$L), M.A. Kaashoek, L. Lerer, \\ and A.L. Sakhnovich}

\date{}
\begin{document}

\maketitle

\begin{abstract}
This paper is devoted to connections between accelerants and potentials of Krein 
systems and of canonical systems of Dirac type,
both on  a finite interval. It is shown that a continuous potential is always generated 
by an accelerant, provided  the latter is continuous 
with a possible jump discontinuity at the origin. Moreover, the generating accelerant is 
uniquely determined by the potential. The 
results are illustrated on pseudo-exponential potentials. The paper is a continuation of the earlier paper of the authors 
\cite{AGKLS07} dealing with the direct problem for Krein systems.
\end{abstract}

\noindent\textit{When Israel Gohberg, the second author of this paper, passed away 
on October 12, 2009, the work on this paper was finished, except for the last section 
of which only the first draft existed. The expression \textup{Z$''$L } after his name is  
used in Hebrew and means ``of blessed memory''. }

\section{Introduction} \label{intro}
\setcounter{equation}{0}

Let $\eT>0$, and let $k$ be a scalar continuous function on
the interval $[-\eT,\eT]$ which is hermitian, that is, $k(-t)=\overline{k(t)}$ on
$-\eT\leq t\leq \eT$. Assume that for each $0< \t\leq \eT$ the corresponding
convolution integral operator $T_\t$ on $L^2(0,\,\t)$,
\begin{equation}
\label{intr} (T_\t f)(t)=f(t)-\int_0^\t
k(t-s)f(s)\,ds, \quad 0\leq t\leq \t,
\end{equation}
is invertible, and let  $\g_\t(t,s)$ be the corresponding
resolvent kernel, i.e.,
\begin{equation}
\label{langoustine}
\g_\tau(t,s)- \int_0^\tau
k(t-\xi)\g_\tau(\xi,s)\,d\xi=k(t-s),\quad 0\le t,s\le
\tau.
\end{equation}
Consider the entire functions
\begin{eqnarray}
\label{defP} {\mathcal P}(\tau,\la)&=& e^{i\la
\tau} \left(1+\int_0^{\tau}e^{-i\la
x}\g_{\tau}(x,0)
dx\right),\\
\label{defP*} {\mathcal
P}_*(\tau,\la)&=&1+\int_0^{\tau}e^{i\la
x}\g_{\tau}(\tau-x, \tau) dx,
\end{eqnarray}
and put $Y(\t,\la)=\begin{bmatrix} {\mathcal
P}(\t,\la)&
{\mathcal P}_*(\t,\la)\end{bmatrix}$. Then,  as was proved by M.G. Krein in
\cite{MR0086978}, the function $Y(\t,\la)$ satisfies
the differential system
\begin{equation}
\label{kr1} \frac{\partial}{\partial \t}
Y(\t,\la)=Y(\t,\la)
\left(i\la\begin{bmatrix}1&0\\0&0
\end{bmatrix}+\begin{bmatrix}0&a(\t)\\ \overline{a(\t)}&
0\end{bmatrix}\right),
\end{equation}
with $a(\t)=\g_\t(\t,0)$ for $\t\in(0,\eT]$. The functions ${\mathcal P}(\t,\la)$ and ${\mathcal P}_*(\t,\la)$ are usually 
referred to as \emph{Krein orthogonal functions}.

We call (\ref{kr1}) a \emph{Krein system}
when, as in the previous paragraph, the function  $a$ is given by  $a(\t)=\g_\t(\t,0)$,
where $\g_\t(t,s)$ is the resolvent kernel corresponding to
some $k$ on $[-\eT, \eT]$ with the properties described in the previous paragraph.
In that case, following Krein,  the function $k$ is called an \emph{accelerant} for  (\ref{kr1}),
and we shall refer to $a$  as the \emph{potential associated with the
accelerant} $k$.

The result referred to above holds in greater generality, namely for systems with
matrix-valued accelerants that are allowed to have a jump
discontinuity at the origin. In fact, in
\cite{AGKLS07} the following result is proved.

\begin{Tm}
\label{oldthm1} Let $k$ be  a ${r\ts r}$-matrix function, which is hermitian, i.e., $k(-t)=k(t)^*$,
and continuous on $-\eT\leq t\leq \eT$ with possibly a jump discontinuity at the origin. Assume that for each $0< \t\leq \eT$ the corresponding integral operator $T_\t$ on $L_r^2(0,\t)$ given by    \eqref{intr} is invertible, and let  $\g_\t(t,s)$ be the corresponding resolvent kernel as
in $(\ref{langoustine})$. Put
\begin{eqnarray}
\label{defP_mat} {\mathcal P}(\tau,\la)&=&
e^{i\la \tau} \left(I_r+\int_0^{\tau}e^{-i\la
x}\g_{\tau}(x,0)
dx\right)\\
\label{defP*_mat} {\mathcal
P}_*(\tau,\la)&=&I_r+\int_0^{\tau}e^{i\la
x}\g_{\tau}(\tau-x, \tau) dx.
\end{eqnarray}
Then $a(\tau)=\g_\tau(0,\tau)$,
with $0<\tau\le \mathbf T$, extends to a continuous function on $[0,\eT]$ and
$Y(\t,\la)=\begin{bmatrix} {\mathcal
P}(\t,\la)&
{\mathcal P}_*(\t,\la)\end{bmatrix}$
satisfies
\begin{equation}
\label{kr2} \frac{\partial}{\partial \t}
Y(\t,\la)=Y(\t,\la)
\left(i\la\begin{bmatrix}I_r&0\\0&0
\end{bmatrix}+\begin{bmatrix}0&a(\t)\\ a(\t)^*&
0\end{bmatrix}\right).
\end{equation}
\end{Tm}
From \cite{gohkol85} we know that the function $\g_\t(t,s)$ is continuous on the triangles on $0\leq s< t\leq \t$ and $0\leq t<s\leq \t$, and that $\g_\t(t,s)$ has a continuous extension on the closures of each of these triangles. Jumps may appear on the diagonal $0\leq s= x\leq \t$. In particular,  the evaluaton of $\g_\t$ at the point $(\t, 0)$, appearing in Theorem \ref{oldthm1},  is well-defined.

As for the scalar case we call \eqref{kr2} a \emph{Krein system} when the
\emph{potential} $a$ is obtained in the way described in Theorem \ref{oldthm1},  and in that case we
say that $k$ is an \emph{accelerant} for \eqref{kr2}.

In this paper we deal, among other things, with the following inverse problem. Consider the system
\eqref{kr2} and assume that the potential $a$ is a ${r\ts r}$-matrix valued function continuous on $[0,\eT]$.
Does it follow that \eqref{kr2} is a Krein system?  In other words, does there exists
a ${r\ts r}$-matrix valued accelerant $k$ on $[-\eT,\eT]$, with possibly a jump
discontinuity at the origin, such that the potential corresponding to $k$ is the
given potential $a$? If we restrict to continuous accelerants, the answer is negative. For instance (see  \cite{AGKLS07}), the potential
\[
a(\tau)=\dfrac{2i}{1+e^{-2i\t}},
\quad\tau\in[0,\, 1],
\] 
does not have a continuous accelerant. However,  we shall prove that for the larger class of
accelerants introduced here, the answer is affirmative.

Krein systems are  closely related to canonical differential systems of Dirac type.
In fact,
if $Y$ is a ${\mathbb C}^{2r\times 2r}$-valued
solution of the system \eqref{kr2} with
potential $a$, then the function
\begin{equation}
\label{hamilcar}
U(\t,\la)=
e^{-i\t\lambda}Y(\t,-2\bar{\la})^*
\end{equation}
is a solution of the canonical system
\begin{equation}
\label{canon}
-ij\frac{\rm
d}{d\t}U(\t,\lambda)=\lambda U(\t,\lambda) +
\begin{bmatrix} 0&v(\t)\\
v(\t)^*&0\end{bmatrix}U(\t,\lambda),
\end{equation}
where
\begin{equation}\label{a2}
j= \begin{bmatrix}
I_{r} & 0 \\ 0 & -I_{r}
\end{bmatrix} \ands  v(\t)=-ia(\t) \quad (0\leq \t\leq \eT).
\end{equation}
It will be convenient to state our main results in terms of a canonical system rather
than a Krein system.

In this paper we show that a continuous matrix-valued potential $v$ is always generated by an accelerant
(provided a jump discontinuity at the origin is allowed) and that an accelerant is uniquely determined by the potential, that is, if  two accelerants generate the same
potential, then they are equal. In fact, we shall prove the following  theorem.

\begin{Tm}\label{mainthm1} Consider the canonical  system \eqref{canon}, and assume that its
potential  $v$ is  continuous on the interval $[0, \, {\bf T}]$. Then, there
is a unique $r \times r$  matrix function
$k$, which is hermitian, i.e.,  $k(-t)=k(t)^*$, and continuous on $-\eT\leq t\leq \eT$ with possibly
a jump discontinuity at the origin,
such that the following holds: for each $0<\t\leq \eT$ the convolution operator
\begin{equation}
\label{intor} (T_\t f)(t)=f(t)-\int_0^\t
k(t-s)f(s)\,ds, \quad 0\leq t\leq \t,
\end{equation}
is invertible on $L^2_r(0, \, \tau)$, and the potential $v$ of  \eqref{canon}  is given by
\begin{equation} \label{s12}
v(\tau)=-i\g_\t(\t,0), \quad 0<\t\leq \eT.
\end{equation}
Here $\g_\t(t,s)$ is the resolvent kernel corresponding to $T_\t$ as in \eqref{langoustine}.
\end{Tm}

In analogy with the theory of Krein systems, a $r\ts r$ matrix function $k$ with the properties described in the above 
theorem will be called an \emph{accelerant} for the canonical system \eqref{canon}.  In this case we also say that the potential
$v$ is \emph{generated} by the accelerant $k$. Using this terminology, Theorem \ref{mainthm1} just tells us that a canonical system
with a continuous potential has a unique accelerant.

Given a continuous matrix-valued potential $v$ we shall also present
a formula for the fundamental solution of the canonical system \eqref{canon} in terms of the accelerant
generating the potential $v$. The result (Theorem \ref{thmdirmain}  in Section \ref{secprdirmain}) can be
viewed as an addition to Theorem \ref{oldthm1}

The statement in Theorem \ref{mainthm1} about the uniqueness of the accelerant is known and has been
proved in \cite{AGKLS07}  using recent results about the continuous analogue of the resultant
for certain entire matrix functions (see Theorem 1.3 in \cite{AGKLS07} for further details). In this paper
we give a new proof of the uniqueness using a formula for the fundamental solution of the canonical system
\eqref{canon} in terms of a given accelerant, which is presented in Theorem \ref{thmdirmain} below.

For the case when the potential  $v$ is bounded,
bounded accelerants $k$ have been constructed
in \cite{SaA3} following the scheme outlined in Section 8.2 of \cite{SaL3} (see also
 \cite{SaA02}). In the present paper, to prove Theorem~\ref{mainthm1},
the approach of \cite{SaA3} is specified and  developed further for the case of
continuous potentials. Also the material related to Theorem \ref{thmdirmain} below is inspired by and builds on results from Sections
3 and 4 in \cite{SaA3}.

The paper consists of six sections including this introduction. In Section \ref{secprdirmain}
we derive a formula for the fundamental solution of the canonical system \eqref{canon} in terms
of its accelerant. The result is used in in Section \ref{secuniqacc} to give  a new proof
of the  uniqueness of the accelerant given the potential as stated in Theorem \ref{mainthm1}.
The next two sections complete the  proof of Theorem \ref{mainthm1}.  Section \ref{semisep}
has an auxiliary character and is interesting in its own right.
We show that a lower triangular  semi-separable integral operator from a certain class
is similar to the operator of integration and that the corresponding similarity operator can be chosen
in such a way that both this similarity operator and its inverse map  functions with
a continuous derivative into functions with a continuous derivative. This result is then
used in Section \ref{secaccel} to prove
Theorem \ref{mainthm1}.  In the final section the main result of Section \ref{secprdirmain} is specified further for 
pseudo-exponential potentials.

\section{The fundamental solution}\label{secprdirmain}
\setcounter{equation}{0}
Throughout this section $k$ is  a $r \times r$  matrix function  on $[-\eT,\, \eT]$, which is hermitian, i.e.
$k(-t)=k(t)^*$, and $k$ is continuous on $-\eT\leq t\leq \eT$ with possibly
a jump discontinuity at the origin. We assume that  $k$ is an  accelerant for the canonical system
\eqref{canon}. The latter means that for each $0<\t\leq \eT$ the convolution operator \eqref{intr}
is invertible on $L^2_r(0, \, \tau)$, and
the potential $v$ is the $r \times r$  continuous matrix function on $[0,\,\eT]$ determined by $k$ via the formula
\begin{equation}
\label{defpot} v(\t)=-i\g_\t(\t,0), \quad 0<\t\leq\eT,
\end{equation}
where $\g_\t(t,s)$is the corresponding resolvent kernel as
in $(\ref{langoustine})$.

We shall derive (explicitly in terms of the accelerant $k$) the fundamental solution $u(x,\la)$ of \eqref{canon} satisfying the initial condition
\begin{equation}
\label{canonincond} u(0, \la)=Q^*,  \quad \mbox{where}\quad Q=\frac{1}{\sqrt{2}}\left[
\begin{array}{cc} I_r &
-I_{r} \\ I_{r} & I_r
\end{array}
\right].
\end{equation}
For this purpose we need the following $r\ts r$ matrix functions:
\begin{equation}\label{defell12}
\ell_1(x, \la)=e^{2i\la x}\left(I_r-2\int_0^x  e^{-2i\la t}k(t)\,dt\right), \quad \ell_2(x, \la)=e^{2i\la x}I_r.
\end{equation}
Both $\ell_1(\cdot, \la)$ and $\ell_2(\cdot, \la)$ are defined on $[0,\, \eT]$. Note that for $0\leq x\leq \eT$ we have
\begin{equation}
\label{diffell12} \frac{d}{dx}\ell_1(x, \la)=2i\la \ell_1(x, \la)-2k(x), \quad  \frac{d}{dx}\ell_2(x, \la)=
2i\la \ell_2(x, \la).
\end{equation}
The next theorem is the main result of this section.

\begin{Tm}
\label{thmdirmain} Assume that the $r\ts r$ matrix function $k$ is an  accelerant for the canonical system
\eqref{canon}, and  let  $\g_\t(t,s)$ be the corresponding resolvent kernel as
in $(\ref{langoustine})$.  For $0\leq \t \leq \eT$, $\la\in \BC$, and  $j=1,2$ put
\begin{eqnarray}
&&\th_j(\t, \la)=\frac{1}{\sqrt{2}}\,e^{-i\t\la} \left\{\ell_j(\t, \la)+\int_0^\t \g_\t(\t,s)\ell_j(s, \la)\,ds\right\} \label{defth12a}\\
&&\om_j(\t, \la)=\frac{1}{\sqrt{2}}\,e^{-i\t\la} \left\{(-1)^j I_r+\int_0^\t \g_\t(0,s)\ell_j(s, \la)\,ds\right\},\label{dirdefom12a}
\end{eqnarray}
where $\ell_1(\cdot, \la)$ and $\ell_2(\cdot, \la)$ are given by \eqref{defell12}.
Then the $2r\ts 2r$ matrix function $u(\t, \la)$ defined by
\begin{equation}
\label{fundsol} u(\t, \la)=\begin{bmatrix}
\th_1(\t, \la)&\th_2(\t, \la)\\
\noalign{\vskip6pt}
\om_1(\t, \la)&\om_2(\t, \la)
\end{bmatrix}, \quad 0\leq \t\leq \eT,
\end{equation}
is the fundamental solution of \eqref{canon} with initial condition \eqref{canonincond}.
\end{Tm}

Using the definition of  $\ell_2(\cdot, \la)$ in the second part of  \eqref{defell12} we see that
\begin{eqnarray}
\th_2(\t,\la)&=&\frac{1}{\sqrt{2}}\,e^{i\t\la}\Big(I_r+\int_0^\t e^{-2i\la s}\g_\t(\t,\t-s)\,ds\Big), \label{formth2}\\
\om_2(\t,\la)&=& \frac{1}{\sqrt{2}}\,e^{-i\t\la}\Big(I_r+\int_0^\t e^{2i\la s}\g_\t(0,s)\,ds\Big).\label{formom2}
\end{eqnarray}
The formula for $\om_2(\t,\la)$ is immediate from the definition and the one for  $\th_2(\t,\la)$ follows using
the following calculation:
\begin{eqnarray*}
\th_2(\t,\la)&=&\frac{1}{\sqrt{2}}\,e^{-i\t\la}\Big(e^{2i\la\t}I_r+\int_0^\t e^{2i\la s}\g_\t(\t,s)\,ds\Big)\\
&=&\frac{1}{\sqrt{2}}\,e^{-i\t\la}\Big(e^{2i\la\t}I_r+\int_0^\t e^{2i\la(\t- s)}\g_\t(\t,\t-s)\,ds\Big)\\
&=&\frac{1}{\sqrt{2}}\,e^{i\t\la}\Big(I_r+\int_0^\t e^{-2i\la s}\g_\t(\t,\t-s)\,ds\Big).
\end{eqnarray*}
The expressions \eqref{formth2} and  \eqref{formom2} show that the functions $\th_2(\t,\la)$ and $\om_2(\t,\la)$ are closely related to the  Krein orthogonal entire matrix functions $\clp(\t,\la)$ and $\clp_*(\t,\la)$ appearing in Theorem \ref{oldthm1}. In
fact we have
\[
\th_2(\t,\la)=\frac{1}{\sqrt{2}}\,e^{-i\t\la}\clp_*(\t, 2\bar{\la})^*, \quad \om_2(\t,\la)=\frac{1}{\sqrt{2}}\,e^{i\t\la}\clp(\t, 2\bar{\la})^*.
\]
Thus Theorem \ref{thmdirmain}  can be seen as an addition to Theorem \ref{oldthm1}.

\medskip
To prove Theorem \ref{thmdirmain} we need some preliminaries.  In the sequel we write $T$ in place of $T_\eT$. The fact that $T_\t$ is invertible for each $0<\t\leq \eT$ is equivalent to $T$ being strictly positive. The latter property implies that
$T$  factorizes as $T=\G\G^*$, where $\G$  is an invertible lower triangular integral operator,
\begin{eqnarray}
&&(\G f)(x)=f(x)+\int_0^x \g_-(x,t)f(s)\, ds, \quad 0\leq x\leq \eT,\label{defGAdir}\\
&&(\G^{-1} f)(x)=f(x)+\int_0^x \g_-^\ts(x,s)f(s)\, ds, \quad 0\leq x\leq \eT,\label{defGAinv}
\end{eqnarray} 
with both  $\g_-(x,s)$ and $\g_-^\ts(x,s)$ being continuous on $0\leq s\leq x\leq \eT$. We shall refer to
$T=\G\G^*$ as the \emph{$LU$-factorization} of $T$. From \cite{gohkol85} we also know that
\begin{equation}\label{gming}
\g_-^\ts(\t,s)=\g_\t(\t,s), \quad 0\leq s\leq \t\leq \eT.
\end{equation}
We shall need the following three lemmas.

\begin{La} \label{kern2} We have
\begin{equation}
\label{pot1} (\G^{-1}k)(\t)=\big(T_\t^{-1}(k|_{[\,0,\,\t]})\big)(\t)=\g_\t(\t,0), \quad 0<  \t\leq \eT.
\end{equation}
\end{La}
\bpr  We first prove the second equality in \eqref{pot1}. Since $\g_\t(t,s)$ is the resolvent kernel corresponding to $T_\t$, we know
(cf., \eqref{langoustine}) that 
\[
\g_\t(t,s)-\int_0^\t k(t-\xi)\g_\t(\xi,s)\,d\xi=k(t-s), \quad 0\leq s\leq t\leq \t.
\]
This equality holds a.e on $0\leq s\leq t\leq \t$. But then, since both $\g_\t(t,s)$ and $k(t-s)$
are continuous on $0\leq s\leq t\leq \t$,  the above equality holds at each point of $0\leq s\leq t\leq \t$.
In particular, at the point $(t,0)$. Thus
\[
\g_\t(t,0)-\int_0^\t k(t-\a)\g_\t(\a,0)\,d\a=k(t), \quad 0\leq t\leq \t.
\]
This shows that $T_\t\g_\t(\cdot,\, 0)=k|_{[\,0,\,\t]}$, and hence
\begin{equation}
\label{identgak}\big(T_\t^{-1}(k|_{[\,0,\,\t]})\big)(x)=\g_\t(x,0), \quad 0\leq x\leq \t.
\end{equation}
For $x=\t$ this yields the second identity in \eqref{pot1}.

Next we prove the first identity in \eqref{pot1}.  Fix $0<  \t\leq \eT$. Since $k$ is continuous on $[0,\,\eT]$, the function $\G^{-1}k$
is continuous on $[0,\,\eT]$.  From the previous part of the proof we know that $T_\t^{-1}(k|_{[\,0,\,\t]})$ is continuous  on $[0,\,\t]$.  Hence for both functions the evaluation at $\t$ is well-defined. Moreover,
\begin{eqnarray*}
(\G^{-1}k)(\t)&=&k(\t)+\int_0^\t \g_-^\ts(\t,s)k(s)\, ds,\\
(T_\t^{-1}(k|_{[\,0,\,\t]})(\t)&=&k(\t)+\int_0^\t \g_\t(\t,s)k(s)\, ds.
\end{eqnarray*}
According to \eqref{gming} we have $\g_-^\ts(\t,s)=\g_\t(\t,s)$ for $0\leq s\leq \t$, which yields the first equality in \eqref{pot1}.
\epr

\begin{La}
\label{contdiff} Let $f$ be a continuously differentiable $\BC^{r\ts m}$-valued function on the interval $[0,\,\eT]$. Then $\G^{-1}f$ is also continuously differentiable on $[0,\,\eT]$ and
\begin{eqnarray}
&&\left(\frac{d}{d\t}\G^{-1}f\right)(\t)-\left(\G^{-1}\frac{d}{d\t}f\right)(\t)=\label{diffGaminvf}\\
&&\hspace{1.5cm} =\g_\t(\t,0)\left(f(0)+\int_0^\t \g_\t(0,s)f(s)\,ds\right), \quad 0<\t\leq \eT.\nonumber
\end{eqnarray}
\end{La}
\bpr Recall that
\[
\left(\G^{-1}f\right)(\t)=f(\t)+\int_0^\t \g_-^\ts(\t,s)f(s)\, ds=f(\t)+\int_0^\t \g_\t(\t,s)f(s)\, ds.
\]
Next, using the generalized Krein-Sobolev identities in formulas (2.10) and (2.11) of \cite{AGKLS07}, we have
\begin{eqnarray*}
\frac{d}{d\t}\int_0^\t \g_\t(\t,s)f(s)\,ds&=&\frac{d}{d\t}\int_0^\t \g_\t(\t,\t-s)f(\t-s)\,ds\\
&=&\g_\t(\t, 0)f(0)+ \textup(\a),
\end{eqnarray*}
where
\begin{eqnarray*}
\textup(\a)&=&\int_0^\t \frac{d}{d\t}\bigg(\g_\t(\t,\t-s)f(\t-s)\bigg)\,ds\\
&=&\int_0^\t \g_\t(\t, 0)\g_\t(0,\t-s)f(\t-s)\,ds +\\
&&\hspace{2cm}+\int_0^\t \g_\t(\t,\t-s)\left(\frac{d}{d\t}f\right)(\t-s)\,ds.
\end{eqnarray*}
We conclude that
\begin{eqnarray*}
\left(\frac{d}{d\t}\G^{-1}f\right)(\t)&=& \left(\frac{d}{d\t}f\right)(\t)+ \int_0^\t \g_\t(\t, s)\left(\frac{d}{d\t}f\right)( s)\,ds+\\
&&\hspace{1cm}+ \g_\t(\t, 0)\left(f(0)+\int_0^\t \g_\t(0,s)f(s)\,ds\right)\\
&=& \left(\G^{-1}\frac{d}{d\t}f\right)(\t)+\\
&&\hspace{1.5cm} +\g_\t(\t,0)\left(f(0)+\int_0^\t \g_\t(0,s)f(s)\,ds\right).
\end{eqnarray*}
Thus $\G^{-1}f$ is continuously differentiable and \eqref{diffGaminvf} holds. \epr

\begin{La}
\label{leminprod} For each $g\in L_r^2(0,\,\eT)$ we have
\begin{equation}
\label{eqinprod} \int_0^\t (\G^{-1}k)(t)^*(\G^{-1}g)(t)\,dt=\int_0^\t \g_\t(0,t)g(t)\,dt, \quad 0\leq \t\leq \eT.
\end{equation}
\end{La}
\bpr The case when $\t=0$ is trivial. Fix $0< \t\leq \eT$.
Let $\Pi_\t$ be the projection of $L_r^2(0,\,\eT)$ onto $L_r^2(0,\,\t)$
defined by $\Pi_\t f=f|_{[\,0,\,\t]}$. Note that $\Pi_\t^*$ is the canonical embedding of
$L_r^2(0,\,\t)$ into $L_r^2(0,\,\eT)$. Put $P_\t=\Pi_\t^*\Pi_\t$. Then $P_\t$ is the orthogonal projection of $L_r^2(0,\,\eT)$ onto the subspace consisting of all functions in $L_r^2(0,\,\eT)$ with support in $[0,\,\t]$. Since $\G$ and $\G^{-1}$ are lower triangular, we have
\begin{equation}
\label{triangcond} \Pi_\t\G=\Pi_\t\G P_\t, \quad \Pi_\t\G^{-1}=\Pi_\t\G^{-1}P_\t.
\end{equation}
From the second identity in \eqref{triangcond}, the definition of $T_\t$, and the factorization
$T_\eT=\G\G^*$,  we
see that
\[
T_\t= \Pi_\t T_\eT\Pi_\t^*=\Pi_\t \G\G^*\Pi_\t^*
=\Pi_\t \G P_\t\G\Pi_\t^*=\Pi_\t \G \Pi_\t^*\Pi_\t\G^*\Pi_\t^*.
\]
Since, by lower triangularity,  $\Pi_\t \G  \Pi_\t^*=(\Pi_\t \G^{-1}\Pi_\t^*)^{-1}$, we obtain
\begin{equation}
\label{Ttauinv} T_\t^{-1}=\Pi_\t \G^{-*}P_\t\G^{-1}\Pi_\t^*.
\end{equation}
If follows that
\begin{equation}
\label{kinvTtau}
T_\t^{-1}(k|_{[\,0,\,\t]})=\Pi_\t \G^{-*}P_\t\G^{-1}\Pi_\t^*(\Pi_\t k)=\Pi_\t \G^{-*}P_\t\G^{-1}P_\t k.
\end{equation}

Now using \eqref{identgak}, the lower triangularity of $\G^{-1}$,  and the above identities one computes that
\begin{eqnarray*}
&&\int_0^\t (\G^{-1}k)(t)^*(\G^{-1}f)(t)\,dt=\int_0^\eT (P_\t\G^{-1}k)(t)^*(P_\t\G^{-1}f)(t)\,dt\\
&&\hspace{1cm} = \int_0^\eT (P_\t\G^{-1}P_\t k)(t)^*(P_\t\G^{-1}P_\t f)(t)\,dt\\
&&\hspace{1cm} = \int_0^\eT (P_\t\G^{-*}P_\t P_\t\G^{-1}P_\t k)(t)^* f(t)\,dt=\int_0^\eT (\Pi_\t^*T_\t^{-1}\Pi_\t k)(t)^*f(t)\,dt\\
&&\hspace{1cm} = \int_0^\t \left(T_\t^{-1}(k|_{[\,0,\,\t]}\right)(t)^*f(t)\,dt=\int_0^\t \g_\t(0,t)f(t)\,dt.
\end{eqnarray*}
This proves \eqref{eqinprod}.\epr

\medskip
Using the identities \eqref{gming} and \eqref{eqinprod} we see that the expressions for $\th_j(\t,\,\la)$ and $\om_j(\t,\,\la)$ in
\eqref{defth12a} and \eqref{dirdefom12a}, respectively, can be rewritten as follows:
\begin{eqnarray}
&&\th_j(\t, \la)=\frac{1}{\sqrt{2}}\,e^{-i\t\la} \left(\G^{-1}\ell_j(\cdot\,,\la)\right)(\t), \label{defth12b}\\
&&\om_j(\t, \la)=\frac{1}{\sqrt{2}}\,e^{-i\t\la}\ts\nonumber \\
&&\hspace{1cm}\ts\left\{(-1)^j I_r+\int_0^\t (\G^{-1}k)(t)^*\left(\G^{-1}\ell_j(\cdot \,,\la)\right)(t)\,dt\right\}.\label{dirdefom12}
\end{eqnarray}
We are now ready to prove Theorem \ref{thmdirmain}.

\medskip\noindent\textbf{Proof of Theorem \ref{thmdirmain}.}
We split  the proof into three parts.  The first part deals with the initial value condition;
the two other parts concern the proof that $u(\t, \la)$ satisfies the differential equation  \eqref{canon}.

\smallskip\noindent\textit{Part 1.} Using the expressions for  $\th_j$ and $\om_j$ in \eqref{defth12b} and \eqref{dirdefom12}, respectively, and the fact that  $\G^{-1}$ is a lower triangular integral operator 
(see \eqref{defGAinv}), we obtain the following identities
\[
\th_j(0, \la)=\frac{1}{\sqrt{2}}\,\ell_j(0\,,\la)=\frac{1}{\sqrt{2}}\,I_r, \quad \om_j(0, \la)=\frac{1}{\sqrt{2}}\,(-1)^jI_r\quad (j=1,2).
\]
Thus
\[
u(0, \la)= \begin{bmatrix}
\frac{1}{\sqrt{2}}\,I_r&\frac{1}{\sqrt{2}}\,I_r\\
\noalign{\vskip6pt}
-\frac{1}{\sqrt{2}}\,I_r&\frac{1}{\sqrt{2}}\,I_r
\end{bmatrix}=\frac{1}{\sqrt{2}}\begin{bmatrix}
I_r&I_r\\
\noalign{\vskip6pt}
-I_r&I_r
\end{bmatrix}=Q^*.
\]
Hence $u(\t,\la)$ has the desired value at $\t=0$.

To complete the proof it suffices to prove the following differential expressions:
\begin{eqnarray}
&&\frac{d}{d\t}\om_j(\t\,,\la)= -i\la \om_j(\t\,,\la)-iv(\t)^*\th_j(\t,\la), \quad j=1,2, \label{diffeqom}\\
\noalign{\vskip6pt}
&&\frac{d}{d\t}\th_j(\t\,,\la)= i\la \th_j(\t\,,\la)+iv(\t)\om_j(\t,\la), \quad j=1,2.\label{diffeqth}
\end{eqnarray}
The first two identities will be proved in the next part and the other two in the final part.

\smallskip\noindent\textit{Part 2.} In this part we prove \eqref{diffeqom}.
Since $k$ is continuous on $[0,\,\eT]$, the fact that the kernel function of the lower triangular integral operator $G^{-1}$ is continuous on
$0\leq s\leq x\leq \eT$ implies that $\G^{-1}k$ is also continuous on $[0,\,\eT]$. Similarly, using the continuity of $\ell_1(\cdot\,, \la)$ and $\ell_2(\cdot\,, \la)$ on $[0,\,\eT]$, we see that $\G^{-1}\ell_1(\cdot\,, \la)$ and $\G^{-1}\ell_2(\cdot\,, \la)$ are continuous on $[0,\,\eT]$. Thus the functions under the integrals in the definitions of $\om_1(\cdot\,,\la)$ and $\om_2(\cdot\,,\la)$ are continuous. This implies that $\om_1(\cdot\,,\la)$ and $\om_2(\cdot\,,\la)$ are continuously differentiable and
\begin{eqnarray*}
&&
\frac{d}{d\t}\om_j(\t\,,\la)= -i\la \om_j(\t\,,\la)+\\
&&\hspace{1cm}+\frac{1}{\sqrt{2}}\,e^{-i\t\la}(\G^{-1}k)(\t)^*\left(\G^{-1}\ell_j(\cdot \,,\la)\right)(\t), \quad j=1,2.
\end{eqnarray*}
Using \eqref{pot1} and \eqref{defpot} we see that $(\G^{-1}k)(\t)=iv(\t)$.
This together with the expression  of $\th_j$ in \eqref{defth12b} shows that \eqref{diffeqom} holds.

\smallskip\noindent\textit{Part 3.}
In this part we prove \eqref{diffeqth}.
First note that
\begin{equation}
\label{step31}\frac{d}{d\t}\,\th_j(\t, \la)=-i\la\th_j(\t, \la)+ \frac{1}{\sqrt{2}}\,e^{-i\t\la} \left(\frac{d}{d\t}\G^{-1}\ell_j(\cdot\,,\la)\right)(\t), \  j=1,2.
\end{equation}
Applying Lemma \ref{contdiff} with $f=\ell_j(\cdot\,,\la)$, $j=1,2$, yields
\begin{eqnarray*}
&&\left(\frac{d}{d\t}\G^{-1}\ell_j(\cdot\,,\la)\right)(\t)=\left(\G^{-1}\frac{d}{d\t}\ell_j(\cdot\,,\la)\right)(\t)+\\
&&\hspace{1cm}+\g_\t(\t,0)\left(\ell_j(0\,,\la)+\int_0^\t \g_\t(0,s)\ell_j(s,\la)\,ds\right), \quad 0\leq \t\leq\eT.
\end{eqnarray*}
Using  the formula for the potential $v$ in \eqref{defpot} and  the identity in \eqref{eqinprod} we obtain
\begin{eqnarray*}
&&\left(\frac{d}{d\t}\G^{-1}\ell_j(\cdot\,,\la)\right)(\t)=\left(\G^{-1}\frac{d}{d\t}\ell_j(\cdot\,,\la)\right)(\t)
+ iv(\t)+\\
&&\hspace{2cm}+iv(\t)\int_0^\t (\G^{-1}k)(t)^*\big(\G^{-1}\ell_j(\cdot\,,\la)\big)(t)\,dt, \quad 0\leq \t \leq\eT.\\
\end{eqnarray*}
 Next we use the identities in \eqref{diffell12} and $(\G^{-1}k)(\t)=iv(\t)$ to show that
 \begin{eqnarray*}
 &&\frac{1}{\sqrt{2}}\,e^{-i\t\la}\left(\G^{-1}\frac{d}{d\t}\ell_1(\cdot\,,\la)\right)(\t)=\\
 \noalign{\vskip4pt}
 &&\hspace{2cm}=2i\la\frac{1}{\sqrt{2}}\,e^{-i\t\la}\left(\G^{-1}\ell_1(\cdot\,,\la)\right)(\t)-
 2\frac{1}{\sqrt{2}}\,e^{-i\t\la}\left(\G^{-1}k\right)(\t),\\
 \noalign{\vskip4pt}
&&\hspace{2cm}=2i\la \th_1(\t, \la)-\sqrt{2}\, e^{-i\t\la}iv(\t),
\end{eqnarray*}
and
\begin{eqnarray*}
\frac{1}{\sqrt{2}}\,e^{-i\t\la}\left(\G^{-1}\frac{d}{d\t}\ell_2(\cdot\,,\la)\right)(\t)&=&
2i\la\frac{1}{\sqrt{2}}\,e^{-i\t\la}\left(\G^{-1}\ell_2(\cdot\,,\la)\right)(\t)\\
\noalign{\vskip4pt}
&=&2i \la\th_2(\t, \la).
\end{eqnarray*}
Returning to \eqref{step31}, first for $j=1$ and next for $j=2$, we obtain
\begin{eqnarray*}
&&\frac{d}{d\t}\,\th_1(\t, \la)=-i\la\th_1(\t, \la)+ 2i\la\th_1(\t, \la) -\sqrt{2}\, e^{-i\t\la}iv(\t)+\\
&&\hspace{.5cm} +iv(\t)\frac{1}{\sqrt{2}}\,e^{-i\t\la}\Big(I_r+\int_0^\t (\G^{-1}k)(t)^*\big(\G^{-1}\ell_1(\cdot\,,\la)\big)(t)\,dt\Big)\\
&&=i\la\th_1(\t, \la)+ iv(\t)\frac{1}{\sqrt{2}}\,e^{-i\t\la}\Big(-2I_r+I_r+\\
&&\hspace{3cm}+\int_0^\t (\G^{-1}k)(t)^*\big(\G^{-1}\ell_j(\cdot\,,\la)\big)(t)\,dt\Big)\\
&&=i\la\th_1(\t, \la)+ iv(\t)\om_1(\t,\la),
\end{eqnarray*}
and
\begin{eqnarray*}
&&\frac{d}{d\t}\,\th_2(\t, \la)=-i\la\th_2(\t, \la)+ 2i\la\th_2(\t, \la)+\\
&&\hspace{1cm}+iv(\t)\frac{1}{\sqrt{2}}\,e^{-i\t\la}\Big(I_r+\int_0^\t (\G^{-1}k)(t)^*\big(\G^{-1}\ell_2(\cdot\,,\la)\big)(t)\,dt\Big)\\
&&=i\la\th_2(\t, \la)+ iv(\t)\om_2(\t,\la),
\end{eqnarray*}
Thus \eqref{diffeqth} is proved.   \epr

\section{Uniqueness of the accelerant}\label{secuniqacc}
\setcounter{equation}{0}
Let $u(x,\la)$ be the fundamental solution of the canonical system \eqref{canon} satisfying the initial condition
\eqref{canonincond}, and put
\begin{equation}
\label{defth1a} \th(x)=\begin{bmatrix} I_r&0 \end{bmatrix}u(x,\,0), \quad 0\leq x \leq \eT.
\end{equation}
Since the potential $v$ of \eqref{canon} is assumed to be continuous, the function $\th$ is continuously
differentiable on $[0,\, \eT]$. With $\th$ we associate a lower triangular semi-separable integral operator $L$
acting on $L_r^2(0,\,\eT)$, namely
\begin{eqnarray}
&&(Lf)(x)=\th(x)J\int_0^x \th(t)^*f(t)\, dt, \quad 0\leq x\leq \eT,\label{defL1}\\
&&J=\begin{bmatrix}
0&I_r\\
I_r&0
\end{bmatrix}.\label{defJ1}
\end{eqnarray}
Note that $L$ depends only on \eqref{canon};
accelerants do not play a role yet.

The operator $L$ will play an important role in the proof of the uniqueness of the accelerant (in the present section), and also
later on in the construction of the accelerant given a continuous potential (in Section \ref{secaccel} below) .

In this section  $k$ is an accelerant for the canonical system \eqref{canon}, and we will show that $k$ is uniquely
determined by the potential $v$.  First we recall that the statement ``$T$ is a convolution operator  with kernel function $k$'' can be expressed in terms of an intertwining relation involving the \emph{operator of integration} which is the operator $A$ on the space $L_r^2(0, \, \eT)$ defined by
\begin{equation}
\label{defopint}
(Af)(x)=\int_0^x f(t)\,dt\quad (0\leq x\leq \eT).
\end{equation}
In fact, using Theorem 1.2 in Chapter 1 of  \cite{SaL2}, we know that
\begin{eqnarray}
&&AT +T A^*= B JB^*,\  \mbox{where $J$ is defined by \eqref{defJ1}},\label{intertwAT} \\
&& B: \BC^{2r}\to L_r^2(0, \, \eT), \quad B y=\frac{1}{\sqrt{2}}\,\ell(\cdot)y\quad (y\in \BC^{2r}).
\label{defB3}
\end{eqnarray}
Here $\ell$ is the $r\ts 2r$ matrix function given by
\begin{equation}\label{dirdefell}
\ell(x)=\begin{bmatrix}
h(x)&I_r
\end{bmatrix} \ \mbox{with \ $h(x)=\displaystyle{I_r- 2\int_0^x k(t)dt}$}, \quad 0\leq x\leq \eT.
\end{equation}
We shall need the following proposition and an additional lemma.

\begin{Pn}\label{propL}
Let $k$ be an accelerant for the canonical system \eqref{canon}, and let $T=\G\G^*$ be the $LU$-factrization of the corresponding
convolution integral operator $T$. Then the operator $L$ defined by \eqref{defL1} is similar to the operator of integration $A$;
in fact, $L=\G^{-1}A\G$.
\end{Pn}
\bpr We first show that $L+L^*=\G^{-1}B JB^*\G^{-*}$. The fact that $k$ is an accelerant for the canonical system \eqref{canon} allows us to use the results of the previous section. Let $\ell_1(x,\la)$ and $\ell_2(x,\la)$ be the $r\ts r$ matrix functions defined by \eqref{defell12}. Note that
$\ell_1(x,0)$ is equal to $h(x)$, where $h$ is the function appearing in \eqref{dirdefell}, and  
$\ell_2(x,0)=I_r$. It follows that 
\[
\ell(x)= \begin{bmatrix} \ell_1(x,0)&\ell_2(x,0) \end{bmatrix}, \quad 0\leq x 
\leq \eT.
\]
But then we see from \eqref{defth1a}, \eqref{fundsol}, and \eqref{defth12b} that
\begin{equation}\label{thsqrell}
\th(x)=\frac{1}{\sqrt{2}}(\G^{-1}\ell)(x), \quad 0\leq x 
\leq \eT.
\end{equation}
Using the definition of $B$ in \eqref{defB3}, the preceding identity yields $\G^{-1}By=
\th(\,\cdot)y$
for each $y\in \BC^{2r}$.  As $L$ is defined by \eqref{defL1}, we obtain $L+L^*=\G^{-1}B JB^*\G^{-*}$.

Next, since  $T=\G\G^*$, the identity in \eqref{intertwAT}  can be rewritten as
\[
A\G\G^* +\G\G^* A^*= B JB^*.
\]
Multiplying the latter identity from the left by  $\G^{-1}$ and from the right by the operator $\G^{-*}$ yields
\[
\G^{-1}A\G +(\G^{-1}A\G)^*= \G^{-1}B JB^*\G^{-*}.
\]
By the result of the first paragraph,  $L+L^*=\G^{-1}B JB^*\G^{-*}$. It follows that
\begin{equation}
\label{LA1} L-\G^{-1}A\G=(\G^{-1}A\G)^*-L^*.
\end{equation}
Note that the operator in the left hand side of \eqref{LA1} is a lower triangular integral operator
of the first kind, while the operator in the right hand side of \eqref{LA1} is an upper triangular operator of the first kind. Hence both sides are equal to the zero operator. Thus $L=\G^{-1}A\G$. \epr

\begin{La}
\label{lemcommA}Let $R$ be an operator on $L_r^2(0,\,\eT)$ commuting with the operator of integration $A$ given by \eqref{defopint}. If, in addition,
$(Ru)(x)=u$   for each $u\in \BC^r$ and each $0\leq x\leq \eT$, then $R$ is the identity operator on $L_r^2(0,\,\eT)$.
\end{La}
\bpr Let $\b$ be the canonical embedding operator from $\BC^r$ into  $L_r^2(0,\,\eT)$, that is, $\b$ is given by $(\b u)(x)=u$
for each $u\in \BC^r$ and $0\leq x\leq \eT$. Then $R\b=\b$. Since $R$ commutes with operator of integration $A$, we have $RA^n\b=A^nR\b=A^n\b$.
Thus $R$ acts as the identity operator on the closed linear hull $\bigvee_{n=0}^\iy \im A^n\b$. By induction one shows
that $\im \b+\im A\b+\cdots+\im A^n\b$ consists of all $\BC^r$-valued polynomials of degree at most $n$. Since the set of
all $\BC^r$-valued polynomials is
dense in $L_r^2(0,\,\eT)$, we conclude that $\bigvee_{n=0}^\iy \im A^n\b$ coincides with
$L_r^2(0,\,\eT)$, and hence $R$ is identity operator on $L_r^2(0,\,\eT)$. \epr

\medskip
\begin{Tm}
\label{thmacc}
The accelerant  is uniquely determined by the potential.
\end{Tm}
\bpr  By specifying \eqref{formth2} for $\la=0$ we see (using
\eqref{gming})  that
\begin{equation}
\frac{1}{2}\sqrt{2}\, (\G^{-1} I_r)(x)=\th_2(x,0)=\begin{bmatrix}I_r&0\end{bmatrix}u(x,\, 0)
\begin{bmatrix}
0\\ I_r
\end{bmatrix},\quad 0\leq x\leq \eT.
\end{equation}
Thus $\G^{-1} I_r$ depends on the potential $v$ only and not on the
particular choice of the accelerant.

Now fix the potential $v$, and let $\tilde{k}$ be another accelerant determining $v$. Thus $\tilde{k}$
is  a hermitian ${r\ts r}$ matrix function on the interval $-\eT\leq t\leq \eT$, which is is continuous on the interval $-\eT\leq t\leq \eT$ with possibly a jump discontinuity at the origin. Furthermore,
the convolution integral operator $\tilde{T}$ defined by
\[(\tilde{T} f)(x)=f(x)-\int_0^\eT
\tilde{k}(x-s)f(s)\,ds, \quad 0\leq x\leq \eT,
\]
is a strictly positive operator on $L_r^2(0,\,\eT)$. Let $\tilde{\G}\tilde{\G}^*$ be the $LU$-factorization of $\tilde{T}$. Then Proposition  \ref{propL}, together with the fact that $L$ depends  on \eqref{canon} only, shows that $\G^{-1}A\G=\tilde{\G}^{-1}A\tilde{\G}$. In other words, the operator
$\tilde{\G}\G^{-1}$ commutes with the operator $A$. The result of the first paragraph
of the proof yields $\G^{-1}I_r= \tilde{\G}^{-1}I_r$.
Thus the operator $\tilde{\G}\G^{-1}$  commutes with $A$ and $(\tilde{\G}\G^{-1}u)(x)=u$ for each $u\in \BC^r$ and $0\leq x\leq \eT$. According to  Lemma \ref{lemcommA} this implies that $\tilde{\G}\G^{-1}$ is the identity operator on $L_r^2(0,\,\eT)$.
Hence $\tilde{\G}=\G$. But then $\tilde{T}=T$, and thus $\tilde{k}=k$. This proves the uniqueness of the accelerant. \epr

\section{Semi-separable triangular operators similar to the operator of integration} \label{semisep}
\setcounter{equation}{0}
Throughout this section $K$ is a semi-separable lower triangular integral operator
on $L^2_r(0, \,\eT)$, that is, the action of $K$ is given by
\begin{equation} \label{a10}
(Kf)(x)=F(x)\int_0^x G(t)f(t)\,dt, \quad f\in L^2_r(0, \,\eT).
\end{equation}
Here  $F(\cdot)$ and $G(\cdot)$ are   matrix functions of sizes $r \times p$ and $p \times r$,
respectively, and their entries are square summable on the interval $[0, \,\eT]$. In fact, we shall
assume that $F(\cdot)$ and $G(\cdot)$ are continuously differentiable on $[0, \,\eT]$ and such that
\begin{equation} \label{n1}
F(x)G(x) = I_r, \quad 0\leq x\leq \eT.
\end{equation}
A simple example of such an operator is the operator of integration A  on $L^2_r(0, \,\eT)$ defined by \eqref{defopint}.

We shall see that any semi-separable lower triangular operator $K$ satisfying the
conditions referred to above is similar to the operator of integration $A$ and with
a similarity operator of a special kind. The precise result is presented in the next proposition.

\begin{Pn}\label{PnSim1} Let $F$ and $G$ be continuously differentiable,
and assume \eqref{n1} holds. Then the operator $K$ defined by \eqref{a10}
is similar to the operator of integration $A$. More precisely, $K=E A
E^{-1}$
where $E$ is a lower triangular operator of the form
\begin{equation} \label{a12}
(E f)(x)=\rho(x) f(x)+\int_0^x e(x,t)f(t)dt, \quad f\in L^2_r(0, \,\eT).
\end{equation}
Here, $e(x,t)$ is a continuous $r\ts r$ matrix function on $0\leq t\leq x\leq \eT$,
which is zero at $t=0$, and
the $r\ts r$ matrix function $\rho$ is given by
\begin{equation}
\label{defrho}\frac{d}{d x}\rho(x) =F^{\prime}(x)G(x)\rho(x), \quad \rho(0)=I_r.
\end{equation}
Moreover, the operators $E^{\pm 1}$ map functions with a continuous derivative
into functions with a continuous derivative.
\end{Pn}

When $F$ and $G$  are boundedly differentiable and continuous derivatives are
replaced by bounded derivatives, the above proposition is  a particular
case of Theorem 1 in \cite{SaL0}. The restriction to continuously differentiable
$F$ and $G$ is the new element here. Since the above proposition plays an essential
role in the proof of our main theorem, we will present a full proof.

\medskip
In order to prove Proposition \ref{PnSim1} we first make some heuristic remarks explaining
the line of reasoning that we will follow. Assume we have an operator $E$ on $L^2_r(0, \,\eT)$
with all the properties described in Proposition \ref{PnSim1}. In particular, $KE=EA$. By rewriting
this identity in terms of the  kernel functions of the integral operators $A$, $K$, $E$, we
get
\[
F(x)G(t)\rho(t)+F(x)\int_t^x G(s)e(s,t)\,ds=\rho(x)+\int_t^x e(x,s)\,ds.\]
Taking $t=0$ and using $e(x,0)=0$ for each $0\leq x\leq \eT$, we conclude that
\begin{equation}\label{EIr}
(EI_r)(x)= \rho(x)+\int_0^x e(x,s)\,ds=F(x)G(0), \quad 0\leq x\leq \eT.
\end{equation}
Here we view $I_r$ as the $r\ts r$ matrix function on $[0,\,\eT]$ which is identically equal to the
$r\ts r$ identity matrix, and $E$ is applied to $I_r$  column wise.

On the other hand, since $K$ and $A$ are Volterra operators (cf., Section 12.9 in \cite{GGK0})
the identity $KE=EA$ implies that
$(I-\la K)^{-1}E=E(I-\la A)^{-1}$ for each $\la\in \BC$. Here, and in the sequel, $I$ denotes the
identity operator on $L^2_r(0, \,\eT)$. As is well known, for each $f$ in $L_r^2(0,\eT)$ we have
\begin{equation}\label{resolA1}
\left((I-\la A)^{-1}f\right)(x)=f(x)+\la\int_0^x
e^{\la(x-t)}f(t)\,dt,\quad 0\leq x\leq \eT.
\end{equation}
With $f(\cdot)=I_r$ this yields
\begin{equation}
\label{resAI} \left((I-\la A)^{-1}I_r\right)(x)=e^{\la x}I_r, \quad 0\leq x\leq \eT.
\end{equation}
It follows (using the identity \eqref{EIr}) that
\begin{eqnarray*}
&&\rho(x)e^{\la x}I_r+\int_0^x e(x,t) e^{\la t}I_r\, dt=\\
&&\hspace{.3cm}=\big(E( e^{\la\cdot}I_r)\big)(x)=\big(E(I-\la A)^{-1}I_r\big)(x)\\
&&\hspace{.3cm}=\big((I-\la K)^{-1}EI_r\big)(x)=\big((I-\la K)^{-1}F(\cdot)G(0)\big)(x), \quad 0\leq x\leq \eT.
\end{eqnarray*}
Hence in order to find the kernel function $e(x,t)$ it is natural to solve
 the equation $(I-\la K)g(\cdot, \la)=F(\cdot)G(0)$ and
to analyze its solution. This will be done in Lemmas \ref{inteq} and \ref{repr} below.

\medskip
We begin with some preparations
According to the general
theory of  semiseparable integral
operators (see Chapter IX in \cite{GGK1}), the inverse of operator $I - \la K$
is given by
\begin{equation} \label{na15}
 \big((I - \la K)^{-1}f\big)(x)=
f(x)+\int_0^x \eta(x,t,\la)  f(t) dt,
\end{equation}
where   and
\begin{eqnarray} \label{na16}&&
  \eta(x,t,\la )=\la F(x)u_1(x,\la)u_1(t,\la)^{-1}G(t), \quad 0 \leq t \leq x \leq {\bf T}, \\
   \label{1b2}
  &&  \frac{d}{d x}\, u_1(x,\la)= \la G(x)F(x)u_1(x,\la), \quad  0 \leq x \leq {\bf T}, \\
 \label{1b3}
&&  u_1(0,\la)=I_r.
\end{eqnarray}
We also need the $r \times r$ matrix function $\wt u_1(x)$ defined by
\begin{equation} \label{1b1}
\frac{d}{d x}\wt  u_1(x)= - G(x)F^{\prime}(x) \wt u_1(x), \quad  0 \leq x \leq {\bf T}, \\
\quad \wt u_1(0)=I_r.
\end{equation}
We are now ready to prove the first lemma.

\begin{La} \label{inteq}
Let $F$  and  $G$ be continuously differentiable, and assume \eqref{n1} holds.
Let $h$ be the $r\times r$ matrix function defined by $h(x)=F(x)G(0)$ on
$0\leq x\leq \eT$, and let $\rho$ be the
$r\times r$ matrix function given by  \eqref{defrho}. Put
\begin{equation} \label{a13}
g(x, \la)=\rho(x)^{-1}\left((I-\la K)^{-1}h\right)(x), \quad 0\leq x\leq \eT,
\end{equation}
where $(I-\la K)^{-1}$ is applied to $h$ columnwise. Then $g$ satisfies the
following integro-differential equation
\begin{equation} \label{a20}
\frac{d}{d x}\,g(x, \la)-\a(x)\int_0^x \b(t) g(t, \la)dt -\la g(x, \la)=0, \quad g(0, \la)=I_r,
\end{equation}
where $\a$ and $\b$ are the continuous functions on  $[0,\,\eT]$ given by
\begin{eqnarray}
\a(x)&=&\rho(x)^{-1}F^{\prime}(x)\wt u_1(x),\quad  0 \leq x \leq \eT, \label{defagr}\\
\noalign{\vskip4pt}
\b(t)&=&-\wt u_1(t)^{-1}\big(G(t)F^{\prime}(t)G(t)+G^{\prime}(t)\big)\rho(t),\quad 0 \leq t \leq \eT.\label{defbgr}
\end{eqnarray}
\end{La}
\bpr
Put $\wt g(x, \la)=\rho(x)g(x, \la)$. Using (\ref{na15})-(\ref{1b3}), (\ref{a13}), and the definition of the matrix
function $h$, we present $\wt g$ in the form
\begin{eqnarray}
&&\wt g(x, \la)=   F(x)G(0)+  \la F(x)u_1(x,\la)\int_0^x u_1(t,\la)^{-1}G(t) F(t)G(0)dt \label{4a1}   \\
&&=F(x)G(0)-F(x)u_1(x,\la)\int_0^x \frac{d}{dt}\Big(u_1(t,\la)^{-1}G(0)\Big)dt\nonumber \\
&& = F(x)G(0)-F(x)u_1(x,\la)\big(u_1(x,\la)^{-1}-I_r\big)G(0)\nonumber\\
&&=F(x)u_1(x,\la)G(0).\nonumber
\end{eqnarray}
It follows that
\begin{equation} \label{4a2}
 g(x, \la)=\rho(x)^{-1}F(x)u_1(x,\la)G(0).
\end{equation}
Clearly  $g$ is differentiable and
\begin{eqnarray} \label{4a3}&&
\frac{d}{d x}\,g(x, \la)=\rho(x)^{-1}\wt g_x(x, \la)-   \rho(x)^{-1}\rho^{\prime}(x)\rho(x)^{-1}   \wt g(x, \la)
\\
\noalign{\vskip4pt}
&&=\rho(x)^{-1} \big\{
\la F(x)G(x)F(x)+F^{\prime}(x) -F^{\prime}(x) G(x)F(x)\big\}u_1(x, \la)G(0)
     \nonumber \\
\noalign{\vskip4pt}
  &&=\la g(x, \la)+\rho(x)^{-1} F^{\prime}(x) \big(I_p-G(x)F(x)\big)u_1(x, \la)G(0). \nonumber
\end{eqnarray}
Here we took into account the identity \eqref{n1}.
From \eqref{1b1} we see that
\[
\frac{d}{dt}  \wt u_1(t)^{-1}=- \wt u_1(t)^{-1}\left(\frac{d}{dt} \wt u_1(t)\right)\wt u_1(t)^{-1}=\wt u_1(t)^{-1}G(t)F^\prime(t).
\]
Hence
\begin{eqnarray*}
&&\frac{d}{dt}  \left(\wt u_1(t)^{-1}\big(I_p-G(t)F(t)\big)
u_1(t,\la)\right)=\\
&&\hspace{1.5cm}  =\wt u_1(t)^{-1}G(t)F^{\prime}(t)\big(I_p-G(t)F(t)\big)u_1(t,\la)+\\
&&\hspace{2.5cm}+\wt u_1(t)^{-1}\big(-G^{\prime}(t)F(t)- G(t)F^{\prime}(t)\big)u_1(t,\la)\\
&&\hspace{2.5cm}+\la\wt u_1(t)^{-1}\big(I_p-G(t)F(t)\big)G(t)F(t)u_1(t,\la).
\end{eqnarray*}
Since $\big(I_p-G(t)F(t)\big)G(t)=0$ because of condition \eqref{n1}, we see that
\begin{eqnarray*}
&&\frac{d}{dt}  \left(\wt u_1(t)^{-1}\big(I_p-G(t)F(t)\big)
u_1(t,\la)\right)=\\
&&\hspace{1.5cm} =\wt u_1(t)^{-1}\big(G(t)F^{\prime}(t)-G(t)F^{\prime}(t)G(t)F(t)-\\
&&\hspace{3.5cm}-G^{\prime}(t)F(t)- G(t)F^{\prime}(t)\big)u_1(t,\la)\\
&&\hspace{1.5cm}=-\wt u_1(t)^{-1}\big(G(t)F^{\prime}(t)G(t)+G^{\prime}(t)\big)F(t)u_1(t,\la).
\end{eqnarray*}
Using the definition of $\b$ in \eqref{defbgr} and the identity \eqref{4a2}, we obtain
\begin{equation}
\frac{d}{dt}  \left(\wt u_1(t)^{-1}\big(I_p-G(t)F(t)\big)
u_1(t,\la)\right)G(0)=\b(t)g(t,\la).
\end{equation}
Recall that  $\big(I_p-G(t)F(t)\big)G(t)=0$ and $\big(I_p-G(0)F(0)\big)G(0)=0$, in particular.
From integration by parts it follows that
\begin{eqnarray*}
&&\int_0^x \b(t)g(t,\la)\,dt=\wt u_1(x)^{-1}\big(I_p-G(x)F(x)\big)u_1(x,\la)G(0)-\\
&&-\big(I_p-G(0)F(0)\big)G(0)-\la \int_0^x\wt u_1(t)^{-1}\big(I_p-G(t)F(t)\big)G(t) \\
&& \times F(t)
u_1(t,\la)G(0)dt= \wt u_1(x)^{-1}\big(I_p-G(x)F(x)\big)u_1(x,\la)G(0).
\end{eqnarray*}
But then, using \eqref{4a3} and the definition of $\a$ in \eqref{defagr}, we arrive at the
identity \eqref{a20}.
\epr

\medskip
The following lemma provides an integral representation of $g$.
\begin{La} \label{repr}
Let $\g(x,t)$ be a $r\ts r$ matrix function, continuous on the interval\  $0\leq t\leq x\leq\eT$.
Then the
integro-differential equation
\begin{equation} \label{b20}
\frac{d}{d x}\,g(x, \la)-\int_0^x \g(x,t) g(t, \la)dt -\la g(x, \la)=0, \quad g(0, \la)=I_r,
\end{equation}
has a unique continuously differentiable
solution $g$. Moreover, $g$ is of the form
\begin{equation}
\label{n3}g(x,\la)=e^{\la x}I_r+\int_0^x e^{\la t} N(x,t) \,dt,\quad
0\leq x\leq \eT,
\end{equation}
where $N(x,t)$ is  continuous on $0\leq t\leq x\leq\eT$ and $N(x,0)=0$ for each $0\leq x\leq \eT$.
\end{La}
\bpr Throughout the proof  we fix $\la\in \BC$. By definition a
solution $g(x,\la)$ of \eqref{b20} is absolutely continuous on
$0\leq x\leq\eT$. In that case, since $\g(x,t)$ is continuous on
$0\leq t\leq x\leq\eT$, we see that
\[
\la g(x,\la)+\displaystyle{\int_0^x\g(x,t)g(t,\la)\,ds} \quad\mbox{is
continuous on $0\leq x\leq \eT$}.
\]
But then $\frac{d}{dx}\,g(x,\la)$ is also continuous on $0\leq x\leq
\eT$. Thus any solution of \eqref{b20} is automatically
continuously differentiable.

By integrating the equation in \eqref{b20} over $0\leq x\leq \t$, where $0\leq \t\leq \eT$, we
obtain the equation
\begin{equation}
\label{b20a} g(\cdot, \la)- \la Ag(\cdot, \la)-AR g(\cdot, \la)=I_r.
\end{equation}
Here $A$ is the operator of integration on $L_r^2(0,\eT)$ defined by \eqref{defopint} and $R$
is the operator on $L_r^2(0,\eT)$ given by
\[
(Rf)(x)=\int_0^x \g(x,t)f(t)\, dt, \quad f\in L_r^2(0,\eT).
\]
Note that for each $f$ in $L_r^2(0,\eT)$ the function $Af$ is
absolutely continuous on $[0,\eT]$. It follows that any solution of
\eqref{b20a} is absolutely continuous, and thus the problems
\eqref{b20} and \eqref{b20a} are equivalent. Using \eqref{resolA1}
we get
\begin{equation}\label{resolA2}
\left((I-\la A)^{-1}Af\right)(x)=\int_0^x e^{\la(x-t)}f(t)\,dt,
\quad 0\leq x\leq \eT.
\end{equation}
From \eqref{resolA2} and the definition of $R$ we see that for each $f$ in $L_r^2(0,\eT)$
\[
\left((I-\la A)^{-1}ARf\right)(x)= \int_0^x \wt{\g}(x,t;\la)f(t)\,dt, \quad 0\leq x\leq \eT,
\]
where
\[
\wt{\g}(x,t;\la)=\int_t^x e^{\la(x-s)}\g(s,t)\,ds, \quad 0\leq t\leq x\leq \eT.
\]
It follows that $I-(I-\la A)^{-1}AR$ is an invertible operator on $L_r^2(0,\eT)$.
Hence the problem \eqref{b20a} has a unique solution in
$L_r^2(0,\eT)$, namely
\begin{equation}
\label{defg1} g(\cdot,\la)=\left(I - (I-\la
A)^{-1}AR\right)^{-1}(I-\la A)^{-1}I_r.
\end{equation}
We conclude that equation \eqref{b20} has a unique continuously differentiable
solution.

It remains to show that the solution $g(\cdot,\la)$ is of the form
\eqref{n3}. To do this, we use \eqref{resAI} and
rewrite  \eqref{defg1}  as
\begin{equation}
\label{defg2}g(\cdot,\la)=e^{\la\cdot}I_r+ \sum_{k=1}^\iy \left((I-\la
A)^{-1}AR\right)^k e^{\la\cdot}I_r.
\end{equation}
Let us compute the first term with $k=1$. Using \eqref{resolA2} and the definition of $R$, we get
\begin{eqnarray*}
&&\left(\big((I-\la A)^{-1}AR\big) e^{\la\cdot}I_r\right)(x)=
\int_0^x e^{\la(x-r)}\left(\int_0^r \g(r,t)e^{\la t}I_rdt\right)  dr=\\
&&\hspace{1cm} =\int_0^x \left(\int_0^r e^{\la(x+t-r)}\g(r,t)\, dt\right)  dr \\
&&\hspace{1cm}=\int_0^x\left(\int_{x-r}^x e^{\la t}\g(r,r+t-x) dt\right)  dr\\
&&\hspace{1cm}=\int_0^x e^{\la t}\left(\int_{x-t}^x\g(r,r+t-x)dr\right)dt= \int_0^x e^{\la t}\g_1(x,t)\, dt,
\end{eqnarray*}
where
\begin{equation}
\label{a24}
\g_1(x,t)=\int_{x-t}^x\g(r,r+t-x)dr, \quad 0\leq t\leq x\leq\eT.
\end{equation}
Next define matrix $\g_k(t,s)$, $k=2,3, \ldots$, recursively by
\begin{equation}
\label{a25}
\g_{k+1}(x,t)=\int_{x-t}^x
\int_{y+t-x}^y \g(y,s)\g_k(s,t+y-x)\,ds\,dy.
\end{equation}
Then, using  similar calculations  as  for $k=1$ above, one proves by induction
that for each $k\geq 1$ we have
\begin{equation}\label{gk23}
\left(\big((I-\la A)^{-1}AR\big)^k
e^{\la\cdot}I_r\right)(x)=\int_0^xe^{\la t}\g_k(x,t)\,dt,\quad 0\leq
x\leq \eT.
\end{equation}
Observe that for each $k$ the function $\g_k(x,t)$ is continuous on $0\leq t \leq x\leq \eT$. Furthermore,
as we see from \eqref{a24} and \eqref{a25}, we have
\begin{equation}\label{gk2}
\|\g_k(x,t)\|\leq c^k\frac{x^{2k-1}}{(2k-1)!}, \quad 0\leq t \leq
x\leq \eT, \quad k\geq 1.
\end{equation}
Here $c$ is a constant independent of $k$. Finally, using \eqref{defg2}, we conclude that \eqref{n3}
holds with
\[
N(x,t)=\sum_{k=1}^\iy \g_k(x,t), \quad 0\leq t\leq x\leq \eT.
\]
By \eqref{gk2} the convergence in the preceding formula is uniform on the triangle $0\leq t\leq x\leq \eT$.
Since each of the terms $\g_k(x,t)$ is continuous on this triangle, it follows that $N(x,t)$
is continuous on  $0\leq t\leq x\leq \eT$ as desired. Finally, from \eqref{a24} and \eqref{a25} it is clear
that $\g_k(x,0)=0$ for each $0\leq x\leq \eT$ and each positive integer $k$. But then $N(\cdot,0)$ is identically equal
to zero too.
\epr

\medskip
\noindent\textbf{Proof of Proposition \ref{PnSim1}}. We split the proof into three parts. In the first part
we define the operator $E$ and establish the similarity $KE=EA$. In the two other parts we prove that
$E^{\pm 1}$ map
functions with a continuous derivative
into functions with a continuous derivative.

\smallskip\noindent\textit{Part 1.} Let $g(x,\la)$ be the matrix function defined by \eqref{a13}.
From Lemma \ref{inteq} we know that $g(x,\la)$ satisfies the integro-differential equation
\eqref{a20}. But then we can apply Lemma \ref{repr} with $\g(x,t)=\a(x)\b(t)$, where $\a(\cdot)$
and $\b(\cdot)$ are defined by \eqref{defagr} and \eqref{defbgr}. It follows that
$g$ admits the representation
\begin{equation}
\label{n3c}
g(x,\la)=e^{\la x}I_r+\int_0^x N(x,t)e^{\la t}I_r\,dt,\quad
0\leq x\leq \eT,
\end{equation}
with $N(x,t)$ being  continuous on $0\leq t\leq x\leq\eT$ and with $N(\cdot,0)$ identically equal to
zero. Now let $E$ be the operator
on  $L^2_r(0, \,\eT)$ defined by
\begin{equation}
\label{a27} (Ef)(x)= \rho(x)f(x)+\int_0^x \rho(x)N(x,t)f(t)\,dt, \quad 0\leq x\leq \eT.
\end{equation}
Here the $r\ts r$ matrix function $\rho$ is defined by \eqref{defrho}. Thus $E$ has the form
\eqref{a12} with $e(x,t)=\rho(x)N(x,t)$. Obviously $e(x,t)$ is continuous on $0\leq t\leq x\leq\eT$
and $e(\cdot,0)=0$  on $[0,\, \eT]$. We claim that this operator $E$ has all the properties described in
Proposition \ref{PnSim1}.

From \eqref{resAI}, \eqref{a13}, formula \eqref{a27} applied to $f=e^{\la\cdot}I_r$, and the identity in \eqref{n3c}
we see that
\[
E(I-\la A)^{-1}I_r=E(e^{\la\cdot}I_r)=\rho(\cdot)g(\cdot,\la)=(I-\la K)^{-1}h,
\]
where $h=F(\cdot)G(0)$. Taking $\la=0$ in the above identity, we obtain $h=EI_r$. Therefore,
\begin{equation} \label{a30}
(I-\la K)^{-1}E I_r=E(I- \la A)^{-1}I_r.
\end{equation}
From the series expansion in (\ref{a30}) it follows that
\begin{equation} \label{a31}
K^j E I_r=E A^{j}I_r, \qquad j=0,1,2, \ldots.
\end{equation}
Therefore, for each $j=0,1,2,\ldots$, we have
 \begin{equation} \label{a32}
(K E)A^{j}I_r=K(EA^{j}I_r)=K^{j+1}E I_r=E A^{j+1}I_r=(E A )A^{j}I_r.
\end{equation}
As the closed linear span  of the columns of the matrices
$\{A^jI_r\}_{j=0}^{\infty}$ coincides
with $L^2_r(0,\, \eT)$, the equalities  in (\ref{a32}) yield   $KE=EA$. Since $E$ is invertible, we obtain
$K=EAE^{-1}$, and hence  $K$ and $A$ are similar. It remains to prove that $E^{\pm 1}$ map
functions with a continuous derivative
into functions with a continuous derivative.

\smallskip\noindent\textit{Part 2.} To show that $E$ has this property, let $f$ be any $\BC^n$-valued function
on $[0,\,\eT]$ with a continuous derivative. Then $f(\cdot)=(Ag)(\cdot)+u$, where $g$ is the derivative of $f$ and $u$
is a constant $r\ts r$ matrix. As we have seen in the previous paragraph, $EI_r=h=F(\cdot)G(0)$. Thus
$Eu=F(\cdot)G(0)u$. According to our hypotheses, $F(\cdot)$ is continuously differentiable.
Hence the same holds true for $Eu$. Next note that
\[
(EAg)(x)=(KEg)(x)=F(x)\int_0^x G(t)(Eg)(t)\,dt.
\]
Since $\rho$ is continuous on $[0,\,\eT]$ and $e(x,t)$ is continuous on $0\leq t\leq x\leq \eT$, we know
that $E$ maps continuous functions into continuous functions. In particular,
$Eg$ is continuous, and hence the above formula shows that $EAg$  has a continuous
derivative. Therefore, $Ef=EAg+Eu$ is continuously differentiable.

\smallskip\noindent\textit{Part 3.} Next, we prove that $E^{-1}$ maps functions with a continuous derivative
into functions with a continuous derivative. First notice that
$E^{-1}$ admits the representation
 \begin{equation} \label{a34}
(E^{-1}f)(x)=\rho(x)^{-1}f(x)  + \int_0^x e^\ts(x,t)f(t) \, dt, \quad f\in L_r^2(0,\,\eT).
\end{equation}
As $e(x,t)$ is continuous on $0\leq t\leq x\leq \eT$, the same holds true for $e^\ts(x,t)$,
and thus $E^{-1}$
maps  continuous function into continuous functions. In terms of the kernel functions the identity
 $E^{-1}E=I$ means
\[
e^\ts(x,t)\rho(t)+\rho(x)^{-1}e(x,t)+\int^x_t e^\ts(x,s)e(s,t)\,ds=0, \quad 0\leq t\leq x\leq \eT.
\]
Recall that $e(\cdot,0)\equiv 0$. Thus by taking $t=0$ in the preceding identity we obtain $e^\ts(x,0)=0$
for $0\leq x\leq \eT$.

 We shall need
the operator $K_1$  on $L^2_r(0,\,\eT)$ defined by
\[
(K_1f)(x)=F^\prime(x)\int_0^xG(t)\,dt, \quad f\in L^2_r(0,\,\eT).
\]
Here $F^\prime$ is the derivative of $F$, which is a continuous function on $[0,\,\eT]$.
Notice that  $K=A(I+K_1)$. Since
$KE=EA$, we have $E^{-1}K=AE^{-1}$ which yields
\begin{eqnarray}
E^{-1}A&=&E^{-1}A(I+K_1)(I+K_1)^{-1}=E^{-1}K(I+K_1)^{-1}\label{a35}\\
&=&AE^{-1}(I+K_1)^{-1}.\nonumber
\end{eqnarray}
Since the kernel function $F^\prime(x)G(t)$ of $K_1$ is continuous,
$(I+K_1)^{-1}$ maps  continuous functions into continuous functions.

Now
let $f$ be a $\BC^r$-valued function
on $[0,\,\eT]$ with a continuous derivative. As in the previous part, we can represent
$f$ as $f(\cdot)=(Ag)(\cdot)+u$, where $g$ is the derivative of $f$ and $u$
is a constant $r\ts r$ matrix. According to \eqref{a35} we have
$E^{-1}Ag=AE^{-1}(I+K_1)^{-1}g$. Since both $E^{-1}$ and $(I+K_1)^{-1}$ map continuous functions
into continuous functions, the function $E^{-1}(I+K_1)^{-1}g$ is continuous. Thus $E^{-1}Ag$ has a continuous
derivative.

Hence in order to prove that $E^{-1}f$ has a continuous
derivative, it suffices to show that $E^{-1}I_r$ has a continuous
derivative. By rewriting the identity $E^{-1}K=AE^{-1}$ in terms of the kernel functions
of $A$, $K$, and $E^{-1}$ we get
\[
\rho(x)^{-1}F(x)G(t)+\int_t^x e^\ts(x,s)F(s)dsG(t)=\rho(t)^{-1}+\int_t^x{e^\ts}(s,t)ds.
\]
By taking $t=0$ and using $e^\ts(x,0)=0$ for $0\leq x\leq \eT$ we obtain
\begin{equation} \label{n8}
I_r-\rho(x)^{-1}F(x)G(0)-\int_0^x{e^\ts}(x,s)F(s)dsG(0)=0.
\end{equation}
Since $F(0)G(0)=I_r$, we can use  (\ref{a34}) and  (\ref{n8}) to show that
\begin{eqnarray} \label{a37}
&&(E^{-1} I_r)(x)=I_r-\rho(x)^{-1}\big(F(x)-F(0)\big)G(0)-
\\ \nonumber
&&\hspace{3cm}-\int_0^x{e^\ts}(x,s)\big(F(s)-F(0)\big)\,ds\,G(0).
\end{eqnarray}
Using $\big(AF^\prime(\cdot)G(0)\big)(x)=\big((F(x)-F(0)\big)G(0)$, it follows
from (\ref{a37}) and  (\ref{a34}) that $(E^{-1} I_r)(x)=I_r-
\big(E^{-1}A F^{\prime}(\cdot)\big)(x)G(0)$.
As  the right-hand side of the latter identity
has a continuous derivative, we obtain that $E^{-1}I_r$ is continuously differentiable.
\epr

\section{Construction of an accelerant}\label{secaccel}
\setcounter{equation}{0}
In this section we establish the main part of  Theorem \ref{mainthm1}. Throughout the $2r\ts 2r$ matrix function $u(x,\la)$
is the fundamental solution of  the  canonical system \eqref{canon} normalized by
\begin{equation}\label{a4}
u(0, \la)=Q^*,  \quad \mbox{where}\quad Q=\frac{1}{\sqrt{2}}\left[
\begin{array}{cc} I_r &
-I_{r} \\ I_{r} & I_r
\end{array}
\right].
\end{equation}
The potential $v$ of \eqref{canon} is assumed to be continuous on $[0,\,\eT]$. Finally, $j$ and $J$  are signature matrices, $j$ is defined by \eqref{a2}  and $J$  by \eqref{defJ1}. Our aim is to show that $v$ is generated by an accelerant.

In what follows $\th$ and $\om$ are the  $r\ts 2r$ matrix functions on $[0,\,\eT]$ defined by
\begin{eqnarray}
\th(x)=\begin{bmatrix}
I_r&0
\end{bmatrix}u(x,0), \quad 0\leq x\leq \eT, \label{defth2a}\\
\noalign{\vskip4pt}
\om(x)=\begin{bmatrix}
0&I_r
\end{bmatrix}u(x,0), \quad 0\leq x\leq \eT.\label{defom2a}
\end{eqnarray}
We begin with two lemmas. The first will enable us to use Proposition~\ref{PnSim1}.

\begin{La}\label{lemth1}
Let $\th$ be the $r\ts 2r$ matrix function on $[0,\,\eT]$ defined by \eqref{defth2a}. Then $\th$
is continuously differentiable on $[0,\,\eT]$,
\begin{equation}\label{propsth}
\th(x)J\th(x)^*=I_r\quad \mbox{and}\quad  \th^\prime(x)J\th(x)^*=0 \quad (0\leq x\leq \eT).
\end{equation}
\end{La}
\bpr
It is straightforward to check that
$Q$ defined in \eqref{a4} satisfies the identities
\begin{equation} \label{a5}
Q^*=Q^{-1}, \quad QjQ^*=J, \quad Q^*JQ=j.
\end{equation}
Since $u(x, \la)$ satisfies   \eqref{canon} and the potential $v$ is continuous, the function
$u(x,\la)$ is continuously differentiable in $x$. Furthermore, again using that
$u(x, \la)$ satisfies   \eqref{canon}, we have $\frac{d}{d x}\left(u(x,
\ov \la)^*ju(x, \la)\right)=0$.
Hence, taking into account (\ref{a4}) and (\ref{a5}), we derive
\begin{equation}   \label{a7}
u(x, \ov \la)^*ju(x, \la)=J, \quad  u(x, \la)Ju(x, \ov
\la)^*=j.
 \end{equation}

The fact that $u(x,\la)$ is continuously differentiable in $x$,
implies  that $\th$ is
continuously differentiable on $[0,\,\eT]$. Furthermore,
from   (\ref{a7})  and \eqref{canon} we see that for each $0\leq x\leq \eT$ we have
\begin{eqnarray}
&&\th(x)J\th(x)^*=\begin{bmatrix}I_r & 0\end{bmatrix}u( x,0)Ju( x,0)^*\begin{bmatrix}
I_r\\0
\end{bmatrix}=I_r, \label{a91} \\
&&\th^{\prime}(x)J\th(x)^*=\begin{bmatrix}I_r & 0\end{bmatrix}
\left(\frac{d}{dx}u( x,0)\right)Ju( x,0)^*\begin{bmatrix}I_r\\0
\end{bmatrix}=0.\label{a92}
\end{eqnarray}
Thus the identities in \eqref{propsth} hold.
\epr

\begin{La}\label{lemom1}
Let $\om$  be the $r\ts 2r$ matrix function defined by \eqref{defom2a}, and let $\th$ be as in
\eqref{defth2a}. Then $\om$ is continuously differentiable on $[0,\,\eT]$, and for each $0\leq x\leq \eT$ we have
\begin{equation}\label{propom3}
\th(x)J\om(x)^*=0, \quad \om^\prime(x)J\om(x)^*=0, \quad \om(0)=\frac{1}{\sqrt{2}} 
\begin{bmatrix}
-I_r&Ir
\end{bmatrix}.
\end{equation}
Moreover, the three identities in \eqref{propom3} determine $\om$ uniquely. Finally,
\begin{equation}\label{omthv}
\om^\prime(x)J\th(x)^*=-iv(x)^*, \quad 0\leq x\leq \eT.
\end{equation}
\end{La}
\bpr Since $u(x,\la)$ is continuously differentiable in $x$ on $[0,\, \eT]$, the same holds true for $\om(x)$.
From the second identity in (\ref{a7}) and the definitions
of $\th$ and $\om$ in \eqref{defth2a} and in \eqref{defom2a}, respectively, we get
\[
\th(x)J\om(x)^* = \begin{bmatrix}
I_r&0
\end{bmatrix}u(x,0)Ju(x,0)^*\begin{bmatrix}
0\\I_r\end{bmatrix}=\begin{bmatrix}
I_r&0
\end{bmatrix}j\begin{bmatrix}
0\\I_r\end{bmatrix}=0.
\]
Analogously, using \eqref{canon},
\begin{eqnarray*}
\om^\prime(x)J\om(x)^*&=&\begin{bmatrix}
0&I_r\end{bmatrix}\left(\frac{d}{dx}\,u(x,0)\right)Ju(x,0)^*\begin{bmatrix}
0\\I_r\end{bmatrix}\\
&=&\begin{bmatrix}
0&I_r\end{bmatrix}\begin{bmatrix}
0&iv(x)\\
-iv(x)^*&0
\end{bmatrix} j \begin{bmatrix}
0\\I_r\end{bmatrix}=0.
\end{eqnarray*}
Thus the first two identities in \eqref{propom3} are proved. 
The third follows directly the normalizing condition 
\eqref{a4}.

To prove that the three identities in  \eqref{propom3} determine $\om$ uniquely, note
that the second identity in \eqref{propom3} implies that
\[
\frac{d}{dx}\left(\om(x)J\om(x)^*\right)=\om^{\prime}(x)J\om(x)^*+\om(x)J\big(\om^{\prime}(x)\big)^*=0.
\]
Thus, using the third identity in  \eqref{propom3}, we obtain $\om(x)J\om(x)^*=-I_r$. 
The latter identity, together with the first identity in  \eqref{propom3}, yields
\begin{equation}
\label{uJom} u(x,0)J\om(x)^*=\begin{bmatrix}0\\-I_r\end{bmatrix}.
\end{equation}
Let $\tilde{\om}$ be another continuously differentiable function on $[0,\,\eT]$ such that 
the three identities in  \eqref{propom3} hold with $\tilde{\om}$ in place of $\om$. 
Repeating the above reasoning for
$\tilde{\om}$ in place of $\om$ we see that \eqref{uJom} holds for $\tilde{\om}$ in place of $\om$.
Thus $u(x,0)J\left(\om(x)^*-\tilde{\om}(x)^*\right)=0$. But the matrices $u(x,0)$ and $J$ 
are non-singular. Thus $\tilde{\om}(x)=\om(x)$ for each $x\in [0,\,\eT]$. 
It follows that  $\om$ is uniquely determined by the identities in \eqref{propom3}.

To prove the final identity \eqref{omthv} we use \eqref{canon} 
and the definitions of $\th$ and $\om $ in \eqref{defth2a} and \eqref{defom2a}. 
This yields
\begin{eqnarray*}
\om^{\prime}(x)J\th(x)^*&=& \begin{bmatrix}
0&I_r
\end{bmatrix}\left(\frac{d}{dx}u(x,0)\right)Ju(x,0)^*\begin{bmatrix}
I_r\\0
\end{bmatrix}\\
&=& \begin{bmatrix}
0&I_r
\end{bmatrix}\begin{bmatrix}
0&iv(x)\\
-iv(x)^*&0
\end{bmatrix}\begin{bmatrix}
I_r\\0
\end{bmatrix}=-iv(x)^*.
\end{eqnarray*}
Hence \eqref{omthv} is proved.\epr

\medskip
In what follows it will be convenient to use the following notation:
\begin{eqnarray}
&&\th_{0, 1}(x)=\theta(x)
\begin{bmatrix}
I_r\\ 0
\end{bmatrix},\quad \th_{0, 2}(x)=\theta(x)
\begin{bmatrix}
0\\ I_r
\end{bmatrix}\quad (0\leq x\leq \eT);\label{thcirc}\\
&&\om_{0,1}(x)=\om(x)
\begin{bmatrix}
I_r\\ 0
\end{bmatrix},\quad \om_{0,2}(x)=\om(x)
\begin{bmatrix}
0\\ I_r
\end{bmatrix}\quad (0\leq x\leq \eT).\label{omcirc}
\end{eqnarray}
Thus
\[
u(x,0)=\begin{bmatrix}\th(x)\\
\noalign{\vskip4pt} \om(x)
\end{bmatrix}=\begin{bmatrix}\th_{0,1}(x)&\th_{0,2}(x)\\ 
\noalign{\vskip4pt}
\om_{0,1}(x)&\om_{0,2}(x)
\end{bmatrix}, \quad 0\leq x\leq \eT.
\]

\medskip
We now return to the operator $L$ defined by \eqref{defL1}.  Thus $L$ is the lower triangular semi-separable
integral operator on $L^2_r(0, \,\eT)$ defined by
\begin{equation}\label{defL2}
(Lf)(x)=\th(x)J\int_0^x \th(t)^*f(t)\, dt, \quad 0\leq x\leq \eT.
\end{equation}
Here $\th$ is as in \eqref{defth2a} (cf., \eqref{defth1a}) and $J$ as in \eqref{defJ1}. Recall (see the 
first paragraph of Section \ref{secuniqacc}) that the definition of $L$ does not involve accelerants 
and depends on \eqref{canon} only. However, if the potential $v$ of  \eqref{canon} is given by an accelerant, then
Proposition \ref{propL} tells us that $L$ is similar to the operator of integration with a similarity operator 
of a special kind. The next proposition goes in the reverse direction. 

\begin{Pn}\label{propLrev} Assume that the operator $L$ defined by \eqref{defL2} 
 is similar to the operator of integration, $L=\Lm^{-1}A\Lm$, where $\Lm$ and $\Lm^{-1}$ have the following
 properties. Both $\Lm$ and $\Lm^{-1}$ are lower triangular operators,
 \begin{eqnarray}
&&(\Lm f)(x)=f(x)+\int_0^x \rho(x,t)f(s)\, ds, \quad 0\leq x\leq \eT,\label{formLm1}\\
&&(\Lm^{-1} f)(x)=f(x)+\int_0^x \rho^\ts(x,s)f(s)\, ds, \quad 0\leq x\leq \eT,\label{formLm2}
\end{eqnarray}  
with  $\rho(x,s)$ and $\rho^\ts(x,s)$ being continuous $r\ts r$ 
matrix functions on the triangles $0\leq s\leq x\leq \eT$.
Furthermore, we assume that $\Lm$ and $\Lm^{-1}$  map continuously differentiable 
functions into continuously differentiable functions, and
\begin{equation}\label{LonI}
(\Lm\th_{0,1})(0)=\frac{1}{\sqrt{2}}Ir, \quad 
\frac{1}{\sqrt{2}}(\Lm^{-1}I_r)(x)=\theta_{0,2}(x)
\quad (0\leq x\leq \eT).
\end{equation}
Then the $r\ts r$ matrix function $k$ given by
\begin{equation}\label{def3k}
k(x)=\left\{
\begin{array}{cl} -\displaystyle{\frac{1}{\sqrt{2}}\frac{d}{dx}}(\Lm \theta_{0,1})(x),& \mbox{for $0< x \leq  \eT$},\\
\noalign{\vskip4pt}
k(-x)^*,&\mbox{for $-\eT\leq x<0$}.
\end{array}\right.
\end{equation}
is an accelerant and $k$ generates  the potential $v$.
\end{Pn}

The above result will allow us to complete the proof of Theorem \ref{mainthm1}. In fact, using 
Proposition \ref{PnSim1}, we shall show that  given $L$ as above  a similarity operator $\Lm$ with 
the properties described in Proposition \ref{propLrev} always exists.

\medskip
\bpr Let $k$ be defined by \eqref{def3k}. Clearly, $k$ is hermitian on $[-\eT,\, \eT]$. 
Since $\theta_{0,1}$ is continuously differentiable, the fact that $\Lm$   maps continuously 
differentiable functions into continuously differentiable functions implies that $k$ 
is continuous on $[-\eT,\, \eT]$ with a possible jump discontinuity at the origin.
The proof that $k$ is an accelerant and generates the potential $v$ will be split into two parts. 

\smallskip
\noindent\textit{Part 1.} In this part we show that $k$ 
is an accelerant. Let $T_\t$ be the operator on $L^2_r(0,\t)$  given
by
\begin{equation}
\label{integr3} (T_\t f)(t)=f(t)-\int_0^\t
k(t-s)f(s)\,ds, \quad 0\leq t\leq \t.
\end{equation}
To prove that $k$ is an accelerant, we have to show that the operator $T_\eT$ 
is strictly positive on $L^2_r(0,\t)$. To establish the 
latter fact we prove the  following identity: 
\begin{equation}
\label{LUfacT3}T= T_\eT =\Lm\Lm^*.
\end{equation}
Note that the right hand side of \eqref{LUfacT3} is an $LU$-factorization.

In order to establish \eqref{LUfacT3}, recall that $L=\Lm^{-1}A\Lm$, 
where $A$ is the operator of integration.
Take $f\in L^2_r(0,\t)$.  Using the similarity relation $L=\Lm^{-1}A\Lm$ it follows that
\begin{equation} \label{a41}
\Lm^{-1} A\Lm f+\big(\Lm^{-1} A\Lm\big)^*f=Lf+L^*f=\th(\cdot)J\int_0^{\eT}\th(t)^* f(t)\, dt.
\end{equation}
By multiplying  (\ref{a41}) from the left by $\Lm$ and replacing $f$ by $\Lm^*f$ we obtain
\begin{equation}\label{a42}
A\Lm \Lm^* f+ \Lm \Lm^*A^*f=(\Lm\th)(\cdot)J\int_0^{\eT}(\Lm \th)(t)^* f(t)\, dt.
\end{equation}
Thus the selfadjoint operator $S=\Lm\Lm^*$, which acts on $L^2_r(0,\t)$, satisfies the identity
\begin{equation}\label{basicid3}
ASf+SA^*f=(\Lm\th)(\cdot)J\int_0^{\eT}(\Lm \th)(t)^* f(t)\, dt,\quad f\in L^2_r(0,\t).
\end{equation}

Now, with $k$ given by \eqref{def3k}, let $s$  be the $r\ts r$  matrix function defined by
\begin{equation} \label{c1}
s(x)=\frac{1}{2}I_r- \int_0^xk(t)dt \qquad 0<x\leq \eT.
\end{equation}
From the first identity in \eqref{LonI} and the definition of $k$ in \eqref{def3k} 
we see that $\Lm \th_{0,1}=\sqrt{2}\,s$. By applying $\Lm$ to both sides of 
the second identity in \eqref{LonI} we obtain $\Lm \th_{0,2}=(\sqrt{2})^{-1}I_r$. 
Summarizing we have
\begin{equation}\label{Lmth1a}
\Lm \th=\frac{1}{\sqrt{2}}\begin{bmatrix}
2s(\cdot)& I_r
\end{bmatrix}.
\end{equation}
Using the later identity  in  the right hand side of \eqref{basicid3} we obtain
\begin{eqnarray*}
&&(\Lm\th)(\cdot)J\int_0^{\eT}(\Lm \th)(t)^* f(t)\, dt=\\
&& \hspace{2cm}=\frac{1}{2}\begin{bmatrix}2s(x)&I_r\end{bmatrix}J\int_0^{\bf T}
\begin{bmatrix}
2s(t)^*\\I_r
\end{bmatrix}f(t)
\, dt\\
&& \hspace{2cm}=\frac{1}{2}\begin{bmatrix}2s(x)&I_r\end{bmatrix}\int_0^{\bf T}
\begin{bmatrix}
I_r\\2s(t)^*
\end{bmatrix}f(t)
\, dt\\
&& \hspace{2cm}=s(x)\int_0^{\bf T}f(t)\,dt+ \int_0^{\bf T}s(t)^*f(t)\,dt, \quad 0\leq x\leq \eT.
\end{eqnarray*}
But then  \eqref{basicid3} can be rewritten as
\begin{equation} \label{a44}
(ASf+SA^*f)(x)= \int_0^{\bf T}\big(s(x)+s(t)^*\big)f(t)\,dt, \quad 0\leq x\leq \eT.
\end{equation}
According to Theorem 2.2 in Chapter 1 of \cite{SaL2} (see also \cite{KKLe94} and \cite{SaLdiff}),
the equation \eqref{a44} has a unique solution which is given by
\begin{equation} \label{a45}
(Sf)(x)=\frac{d}{dx}\int_0^{\bf T}s(x-t) f(t)  \, dt,
\quad s(-x)=-s(x)^* \quad (0< x\leq \eT).
\end{equation}
(Note Theorem 2.2 in Chapter 1 of \cite{SaL2} is stated for scalar kernel functions, but the result also holds for
matrix-valued kernel functions \cite{SaL80b}.  In fact, to get the result for
matrix-valued kernel functions one just writes $S$  as a $r \times r$ matrix  with operator
entries and applies the scalar-valued result to each of these entries.)
From \eqref{a45} and  \eqref{c1} we see that $S=T_\eT$, and thus \eqref{LUfacT3} is proved. In particular,
$k$ is an accelerant.

\smallskip
\noindent\textit{Part 2.}  Let $\tilde{v}$ be the potential generated by the accelerant $k$, where $k$ is as in the
previous part. In this part we show that $v=\tilde{v}$. 

Consider the canonical system \eqref{canon} with the potential $v$ 
being replaced by $\tilde{v}$. 
Let $\tilde{u}(x,\la)$ be the corresponding fundamental solution normalized 
at $x=0$ by $\tilde{u}(0,\la)=Q^*$, where $Q$ is as in \eqref{canonincond}. Put
\[
\tilde{\th}(x)=\begin{bmatrix}
I_r&0
\end{bmatrix}\tilde{u}(x,0), \quad 
\tilde{\om}(x)=\begin{bmatrix}
0&I_r
\end{bmatrix}\tilde{u}(x,0 \quad (0\leq x\leq \eT).
\]
From \eqref{omthv} we know that
\begin{equation}\label{tildev}
\tilde{\om}^\prime J\tilde{\th}^*=-i\tilde{v}^*, \quad 0\leq x\leq \eT.
\end{equation}
Thus to prove $v=\tilde{v}$ it suffices to show that $\th=\tilde{\th}$ and $\om=\tilde{\om}$.

We first show that $\th=\tilde{\th}$. 
Since $k$ is an accelerant generating the potential $\tilde{v}$, we can apply the 
results of Sections \ref{secprdirmain} 
and \ref{secuniqacc}
to the canonical system \eqref{canon} with $\tilde{v}$ in place of $v$. 
In particular, using \eqref{thsqrell}  in the present setting, we see that
\[
\tilde{\th}=\frac{1}{\sqrt{2}}\Lm^{-1}\tilde{\ell}, \quad \mbox{where}\quad  \tilde{\ell}(x)=
\begin{bmatrix}
I_r-2\displaystyle{\int_0^x {k}(t)\,dt}&I_r
\end{bmatrix}.
\]
Here $\Lm$ is the lower triangular integral operator appearing in the 
$LU$-facto--rization \eqref{LUfacT3} of the convolution operator 
$T=T_\eT$ defined by $k$ via \eqref{integr3}. 
By  \eqref{c1} we have $\tilde{\ell}(x)=
\begin{bmatrix}
2s(x)&I_r
\end{bmatrix}$, and hence, using \eqref{Lmth1a}, we obtain  $\tilde{\th}=\th$.

Next we prove that $\tilde{\om}=\om$. By applying Lemma \ref{lemom1} 
to the canonical system \eqref{canon} with $\tilde{v}$ in place
of $v$, we have 
\begin{equation*}
\tilde{\th}(x)J\tilde{\om}(x)^*=0, \quad 
\tilde{\om}^\prime(x)J\tilde{\om}(x)^*=0, \quad \tilde{\om}(0)=\frac{1}{2}\sqrt{2}\,
\begin{bmatrix}
-I_r&Ir
\end{bmatrix}.
\end{equation*}
However,  $\tilde{\th}=\th$. Thus \eqref{propom3} holds with $\tilde{\om}$ in 
place of $\om$. But then we can use the uniqueness
 statement in Lemma  \ref{lemom1} to show that $\tilde{\om}=\om$. 
 
 We have now proved that $v=\tilde{v}$, and hence $k$ is an accelerant generating the potential $v$. \epr

\medskip
\noindent\textbf{Completing the proof of Theorem \ref{mainthm1}.}
Let $L$ be the lower triangular semi-separable integral operator defined by \eqref{defL1}; see also \eqref{defL2}. In order to
complete the proof of Theorem \ref{mainthm1} it suffices to show that $L$ is similar to the operator $A$ of integration, $L=\Lm^{-1}A\Lm$,
where $\Lm$ has all the properties stated in Proposition \ref{propLrev}. For this purpose we use Proposition \ref{PnSim1} with
\begin{equation}\label{KLFG}
K=L, \quad F(x)=\th(x), \quad G(x)=J\th(x)^* \quad (0\leq x \leq \eT).
\end{equation}
By Lemma \ref{lemth1}, the functions $F$ and $G$ in \eqref{KLFG} are continuously differentiable on $[0,\,\eT]$
and condition \eqref{n1} is satisfied. Furthermore, the second identity in \eqref{propsth} implies that 
for $F$ and $G$ in \eqref{KLFG} the solution $\rho$ of the differential  equation \eqref{defrho} is identically equal to $I_r$.
Thus, by Proposition~\ref{PnSim1}, 
\begin{equation} \label{nad}
L=EAE^{-1},
\end{equation}
 where $A$ is the operator of integration defined by \eqref{defopint}
and $E$ on $L_r^2(0,\,\eT)$ is a lower triangular integral operator of the form
\begin{equation} \label{ca12}
(E f)(x)=f(x)+\int_0^x e(x,t)f(t)dt, \quad f\in L^2_r(0, \,\eT).
\end{equation}
Moreover, we know that $e(x,t)$ is a continuous $r\ts r$ matrix function on $0\leq t\leq x\leq \eT$,
which is zero at $t=0$, and the operators $E^{\pm 1}$ map functions with a continuous derivative
into functions with a continuous derivative.

\medskip
To construct the lower triangular integral operator $\Lm$ we  need (apart from the operator $E$) an additional
normalizing lower triangular operator. This operator is the lower triangular convolution operator  $E_0$ defined by 
\begin{eqnarray} \label{a38}
&&(E_0f)(x)=\th_{0,2}(0)f (x) + \int_0^x e_0(x-t) f(t) dt,\quad \mbox{where}\label{a38a} \\
&&\hspace{2cm} e_0(x):=\frac{d}{dx}(E^{-1}\th_{0,2})(x).\label{a38b}
\end{eqnarray}
Recall that $\th_{0,2}$ is defined by the second identity in \eqref{thcirc}. Since $\th$ is continuously differentiable
(see Lemma \ref{lemth1}), the same holds true for $\th_{0,2}$. Using the fact that $E^{-1}$ maps functions with 
a continuous derivative into functions with a continuous derivative, we conclude that $e_0$ is continuous (in fact, continuously
 differentiable).

\begin{La}
\label{lemE0} Let $E_0$ be the operator on $L^2_r(0,\, \eT)$ defined by \eqref{a38}, and let $A$
be the operator of integration defined by \eqref{defopint}. Then
\begin{equation} \label{a39}
E_0A=AE_0 \ands (E_0 I_r)(x)=(E^{-1}\th_{0,2})(x)\quad (0\leq x\leq \eT).
\end{equation}
Furthermore, $E_0$ is invertible and  $E_0^{\pm 1}$ map functions  with a continuous derivative
into functions with a continuous derivative.
\end{La}
\bpr
Since $E_0$ is a lower triangular convolution integral operator, $E_0$
commutes with the operator of integration. Thus the first identity in
\eqref{a39} holds. From \eqref{ca12} with $f=\th_{0,2}$ we see that $(E^{-1}\th_{0,2})(0)=\th_{0,2}(0)$.
By using the latter identity,  \eqref{a38a}, and  \eqref{a38b} we obtain
\begin{eqnarray*}
(E_0 I_r)(x)&=&\th_{0,2}(0)+\int_0^x e_0(x-t)\,dt=\th_{0,2}(0)+\int_0^x e_0(t)\,dt\\
&=&\th_{0,2}(0)+\int_0^x \frac{d}{dt}(E^{-1}\th_{0,2})(t)\, dt\\
&=&\th_{0,2}(0)+(E^{-1}\th_{0,2})(x)-\th_{0,2}(0) =  (E^{-1}\th_{0,2})(x),
\end{eqnarray*}
which yields the second identity in \eqref{a39}.
According to  \eqref{defth2a}, \eqref{thcirc}, and the initial condition in \eqref{a4}, we have
\begin{equation}\label{c2}
\begin{bmatrix}
\th_{0,1}(0)&\th_{0,2}(0)
\end{bmatrix}=\begin{bmatrix}
I_r&0
\end{bmatrix}u(0,0)=\frac{1}{\sqrt{2}}\begin{bmatrix}
I_r & I_r
\end{bmatrix}.
\end{equation}
 In particular, $\th_{0,2}(0)=I_r/\sqrt{2}$,  and so $E_0$ is invertible. Furthermore, $E_0^{-1}$ is
of the form
\begin{equation}
\label{E0x}
(E_0^{-1}f)(x)=\th_2(0)^{-1}f(x)+\int_0^x e_0^\ts(x-t)f(t)\, dt, \quad 0\leq x\leq \eT,
\end{equation}
with $e_0^\ts(x)$ being continuous on $0\leq x\leq \eT$.

Next, let $f$ be any  $\BC^r$-valued function
on $[0,\,\eT]$ with a continuous derivative. Write
$f$ as $f(\cdot)=(Ag)(\cdot)+u$, where $g$ is the derivative of $f$ and $u$
is a constant $r\ts r$ matrix. Then $E_0f=E_0Ag +E_0u=AE_0g +E_0u$. Since $e_0$ and $g$ are continuous
functions, $E_0g$ is continuous, and thus $E_0Ag$ is continuously differentiable. Hence in order to prove
that $E_0f$ is continuously differentiable, it suffices to show that  $E_0u$ has this property.
The latter can   be derived  from the second identity in \eqref{a39} and the
properties of $E$. A more direct argument is as follows. From \eqref{a38} we see that
\[
(E_0u)(x)=\th_2(0)u+\int_0^xe_0(x-t)u\,dt=\th_2(0)u+\int_0^x e_0(t)u\,dt.
\]
Since $e_0$ is continuous, this implies that $E_0u$  is continuously differentiable as desired.
In a similar way, using that $E_0^{-1}$ commutes with $A$ and that $E_0^{-1}$ is given by \eqref{E0x} with
$e_0^\ts$ being continuous, one shows that $E_0^{-1}$ maps functions with a continuous derivative into
functions with a continuous derivative.
\epr

\medskip
For latter purposes we note that
\begin{equation}
\label{atzero} (E_0^{-1}E^{-1}\th_{0,1})(0)=I_r.
\end{equation}
To see this, observe that by \eqref{E0x} for any continuous $\BC^r$-valued function 
$f$ we have $(E_0^{-1}f)(0)=\th_{0,2}(0)^{-1}f(0)$. 
We apply this identity to $f=E^{-1}\th_{0,1}$. We know that $\th_{0,1}$ is continuously differentiable,
and hence $E^{-1}\th_{0,1}$ has the same property. In particular, $E^{-1}\th_{0,1}$
 is continuous. Using \eqref{a34}
and the fact that in this case $\rho$ defined by \eqref{defrho} is identically equal to $I_r$, we see that
$(E^{-1}\th_{0,1})(0)=\th_{0,1}(0)$. Thus
\[
(E_0^{-1}E^{-1}\th_{0,1})(0)=\th_{0,2}(0)^{-1}(E^{-1}\th_{0,1})(0)=\th_{0,2}(0)^{-1}\th_{0,1}(0).
\]
But then \eqref{c2} yields \eqref{atzero}.

Now define 
\begin{equation}
\label{defintopLm} \Lm=\frac{1}{\sqrt{2}} E_0^{-1}E^{-1}.
\end{equation}
We claim that $\Lm$ given by \eqref{defintopLm} satisfies all the conditions on $\Lm$ stated in Proposition \ref{propLrev}.
Indeed, from \eqref{nad} and the first identity in \eqref{a39} we see that $L=\Lm A\Lm^{-1}$. Furthermore, $\Lm$ and $\Lm^{-1}$
are lower triangular integral operators of the form \eqref{formLm1} and \eqref{formLm2}, respectively, and their respective kernel 
functions are continuous on the triangles $0\leq s \leq t\leq \eT$,  because the kernel functions of $E^{\pm 1}$ and $E_0^{\pm 1}$ 
have these properties.  Since $E^{\pm 1}$ and $E_0^{\pm 1}$ map functions with a continuous derivative into functions 
with a continuous derivative, the same holds true for  $\Lm$ and $\Lm^{-1}$. It remains to check the identities in
\eqref{LonI}.  The first identity in \eqref{LonI} follows from the definition of $\Lm$ in \eqref{defintopLm} and the equality in \eqref{atzero}. 
Finally, we use the second equality in \eqref{a39}. The latter can be rewritten as 
$E_0^{-1}E^{-1}\th_{0,2}=I_r$. Using definition of $\Lm$ in \eqref{defintopLm}, this yields  the second identity in \eqref{LonI}.

Thus $\Lm$ given by \eqref{defintopLm} satisfies all the conditions on $\Lm$ appearing  in Proposition \ref{propLrev}. Hence the 
potential $v$ is generated by an accelerant, as desired.\epr

\section{Pseudo-exponential potentials} \label{Exa}
\setcounter{equation}{0}
In this section we consider the class of so-called
pseudo-exponential potentials, which has been introduced in \cite{GKSakh98}; see also \cite{GKSakh02}. 
The aim is to show how Theorem \ref{thmdirmain} can be used to 
present an alternative proof of the basic formula for the fundamental solution given in
Theorem 4.2 of \cite{GKSakh98}; see also Section 2 in \cite{GKSakh02}.

 We begin with some notation. Fix an integer $n>0$ and a triple
of parameter matrices: an $n \times n$ matrix $\clb$ and $n \times r$ matrices $\Phi_1$ and
$\Phi_2$. Recall that the triple  $\clb, \, \Phi_1, \, \Phi_2$ is called {\it admissible} whenever
\begin{equation} \label{e1}
\clb^* -\clb=i \Phi_2\Phi_2^*. 
\end{equation}
Throughout $\Phi$ is the $n \times r$  matrix given by $\Phi=\Phi_1+i\Phi_2$. 

Now let  $\clb$, $\Phi_1$ and $\Phi_2$ be an admissible triple, and put
\begin{equation}
\label{e3} k(t)=-2\Phi_1^*e^{2it\clb^*}\Phi,
\qquad k(-t)=k(t)^*,  \qquad  t>0.
\end{equation}
By taking adjoints, a minor modification  of the proof of Proposition 5.2  in \cite{AGKLS07} 
shows that the function $k$ is an accelerant 
on each interval $[-\eT,\eT]$, and  the corresponding potential is given by
\begin{equation}
\label{e4} v(\t)=2i\Phi_1^*e^{i\t \cla^*}\Sigma(\t)^{-1}e^{i\t \cla}\Phi, \quad \cla=\clb-\Phi_1\Phi_2^*, 
\end{equation}
where
\begin{equation}
\label{e5} \Sigma(t)=I_n+\int_0^t\Pi(s)\Pi(s)^*\, ds, \quad \mbox{with}\quad \Pi(t)=\begin{bmatrix}
e^{-it\cla}\Phi_1&- e^{it\cla}\Phi
\end{bmatrix}.
\end{equation}
Note that with $\Phi_1=\g_1$ and $\Phi_2=\g_2$, we have $\Phi=\g_1+i\g_2$, and in this case  $v$ in \eqref{e4} 
is just equal to $v$ given by (4.6) in \cite{GKSakh98}. The following result is a variant of  Theorem 4.2  in \cite{GKSakh98}. 

\begin{Pn}\label{propfundsol6} Let $\clb$, $\Phi_1$ and $\Phi_2$ be an admissible triple, 
and let $v$ be the potential defined by \eqref{e4}. 
Then the fundamental solution $u(x, \la)$ of the canonical system \eqref{canon}
satisfying the initial condition \eqref{canonincond} is given by 
\begin{equation}\label{e2}
u(\t, \lambda )=w_{ \cla ,
\Pi }(\t, \lambda )e^{i\t \lambda j}w_{ \cla ,
\Pi }(0, \lambda )^{-1}Q^*,
\end{equation}
where $j$ is the $2r \ts 2r$ matrix in the left hand side of \eqref{a2} and 
\begin{equation}\label{e6}
w_{ \cla , \Pi }(\t, \lambda )=I_{2r}+ij \Pi (\t)^{*} \Sigma (\t)^{-1}
( \lambda I_{n}- \cla )^{-1} \Pi (\t).
\end{equation}
\end{Pn}

The proof of Proposition \ref{propfundsol6} given below  is very different from the 
proof of Theorem 4.2 in \cite{AGKLS07}. Here we shall use that the potential 
$v$ in \eqref{e4} is generated by the accelerant $k$ in \eqref{e3}. 
This fact will  allow us to employ the formula for the fundamental solution given in  
Theorem \ref{thmdirmain}. 

We shall only prove equality  \eqref{e2}  for the block $\om_2(\t, \la)$ of $u(\t, \la)$
(see \eqref{fundsol}); the representation of the other blocks can be proved
in a similar way. 

\smallskip
\bpr We shall show that  $\om_2=\wh \om_2$, where  $\wh \om_2$ denotes the right lower block 
on the right-hand side of \eqref{e2}.
In Theorem \ref{thmdirmain}  the block  $\om_2(\t, \la$ is given (cf., \eqref{formom2}) by 
 \begin{equation}\label{e16}
 \om_2(\t, \la)=\frac{1}{\sqrt{2}}\,e^{-i\t\la} \left\{ I_r+\int_0^\t e^{2i s \la }\g_\t(0,s)ds\right\}.
\end{equation}
Here $\g_\t(t,s)$ is the resolvent kernel corresponding to the accelerant $k$. Using adjoints, the same line of reasoning
as in the proof of Proposition 5.2 in \cite{AGKLS07},  shows that
\begin{equation}\label{e17}
\g_{\t}(0,s)=-2\Phi^*e^{-i\t \cla^*} \S(\t)^{-1}e^{-i\t \cla}
\begin{bmatrix}I_n & 0
\end{bmatrix}
e^{i (s-\t) A^{\times}_M}\begin{bmatrix}\Phi_1 \\ \Phi_2
\end{bmatrix},
\end{equation}
where
\begin{equation}\label{e11}
 A^{\times}_M=-2
\begin{bmatrix}
\cla & -\Phi_1\Phi_1^* \\ 0 & \cla^*
\end{bmatrix}.
\end{equation}

In what follows we shall use the identity
\begin{equation}\label{e13}
\begin{bmatrix}I_n & 0
\end{bmatrix}e^{-i \t A^{\times}_M}
\begin{bmatrix}
I_n \\ i I_n
\end{bmatrix}=e^{i\t \cla} \S(\t)e^{i\t \cla^*} .
\end{equation}
Here  $A^{\times}_M$  and $\S(\t)$ are as in \eqref{e11} and  \eqref{e5}, respectively. Note that
\eqref{e13} is the analogue of formula (4.7) in \cite{GKSakh98}.

By substituting  \eqref{e17} in  \eqref{e16} we get
 \begin{eqnarray}\nonumber
 &&\om_2(\t, \la)=\frac{1}{\sqrt{2}}\,e^{-i\t\la} \Big\{ I_r+2i
\Phi^*e^{-i\t \cla^*} \S(\t)^{-1}e^{-i\t \cla}
\begin{bmatrix}I_n & 0
\end{bmatrix}
\\ &&\hspace{.5cm} \times
\big(2\la I_{2n}+A^{\times}_M\big)^{-1}
\Big(e^{2i\t \la }I_{2n}-e^{-i \t A^{\times}_M}\Big)\begin{bmatrix}\Phi_1 \\ \Phi_2
\end{bmatrix}
\Big\}.\label{e18}
\end{eqnarray}
Since $A^{\times}_M$ is given by  \eqref{e11}, we can rewrite \eqref{e18} in the form
 \begin{eqnarray}\nonumber
 &&\om_2(\t, \la)=\frac{1}{\sqrt{2}}\,e^{-i\t\la} \Big\{ I_r-i
\Phi^*e^{-i\t \cla^*} \S(\t)^{-1}e^{-i\t \cla}(\cla -\la I_n)^{-1}
\\ &&\hspace{.5cm} \times
\begin{bmatrix}I_n & \Phi_1\Phi_1^*(\cla^* -\la I_n)^{-1}
\end{bmatrix}
\Big(e^{2i\t \la }I_{2n}-e^{-i \t A^{\times}_M}\Big)\begin{bmatrix}\Phi_1 \\ \Phi_2
\end{bmatrix}
\Big\}.\label{e19}
\end{eqnarray}
Partition $e^{-i \t A^{\times}_M}$ into $n \times n $ blocks $\Big(e^{-i \t A^{\times}_M}\Big)_{kj}$.
From  \eqref{e11} and \eqref{e13} it follows that
\begin{eqnarray}\label{e20}
&& \Big(e^{-i \t A^{\times}_M}\Big)_{11}= e^{2i\t \cla}, \quad \Big(e^{-i \t A^{\times}_M}\Big)_{22}= e^{2i\t \cla^*}, \\
\label{e21} &&
\Big(e^{-i \t A^{\times}_M}\Big)_{21}= 0, \quad \Big(e^{-i \t A^{\times}_M}\Big)_{12}= 
i\Big(e^{2i\t \cla}-e^{i\t \cla} \S(\t)e^{i\t \cla^*}\Big) .
\end{eqnarray}
Taking into account  \eqref{e19}-\eqref{e21} we arrive at
\begin{eqnarray}\nonumber
 &&\om_2(\t, \la)=\frac{1}{\sqrt{2}}e^{i\t\la}\Big\{-i\Phi^*e^{-i\t \cla^*} \S(\t)^{-1}e^{-i\t \cla}(\cla -\la I_n)^{-1}
\Phi_1
\\ \label{e22} && \times
\Big(I_r+\Phi_1^*(\cla^* -\la I_n)^{-1}\Phi_2\Big)\Big\}
+\frac{1}{\sqrt{2}}e^{-i\t\la}
\\ \nonumber && \times
\Big\{I_r+i
\Phi^*e^{-i\t \cla^*} \S(\t)^{-1}e^{-i\t \cla}(\cla -\la I_n)^{-1}
\Big(e^{2i\t \cla}\Phi_1
\\ \nonumber && 
+\Phi_1\Phi_1^*(\cla^* -\la I_n)^{-1}e^{2i\t \cla^*}\Phi_2+i\big(e^{2i\t \cla}-e^{i\t \cla} \S(\t)e^{i\t \cla^*}\big)\Phi_2\Big)
\Big\}.
\end{eqnarray}

Now, consider  the right lower block $\wh \om_2$ of the right-hand side
of \eqref{e2}. The transfer matrix function $w_{ \cla , \Pi }(\t, \lambda )$ 
has the property (see, e.g., \cite{SaL3}):
\[
w_{ \cla ,
\Pi }(\t, \ov \lambda )^*jw_{ \cla ,
\Pi }(\t, \lambda )=j.
\] 
In particular, we have $w_{ \cla ,
\Pi }(0, \lambda)^{-1}=jw_{ \cla ,
\Pi }(0, \ov \lambda )^*j$. Hence, using  \eqref{canonincond},  \eqref{e5}, and  \eqref{e6}, we can write 
\begin{eqnarray}
 &&\wh \om_2(\t, \la)=\nonumber\\
 &&=\frac{1}{\sqrt{2}}
\Big( \begin{bmatrix}0 & I_r
\end{bmatrix}-i\Phi^*
e^{-i\t \cla^*} \S(\t)^{-1}(\cla -\la I_n)^{-1}\begin{bmatrix}
e^{-i\t\cla}\Phi_1&-e^{i\t\cla}\Phi
\end{bmatrix}\Big) \nonumber
\\  && \hspace{1.5cm}\times e^{i\t\la j}\Big(I_{2r}+ij
\begin{bmatrix}
\Phi_1^* \\ -\Phi^*
\end{bmatrix}(\cla^* -\la I_n)^{-1}\begin{bmatrix}\Phi_1& -\Phi
\end{bmatrix}\Big)\begin{bmatrix}
I_r \\ I_r
\end{bmatrix}. \label{e23}
\end{eqnarray}
Formula \eqref{e22} has the form
\begin{equation}\label{e24}
\om_2 (\t,\la)=\frac{1}{\sqrt{2}}\,e^{i\t\la}c_+(\t,\la)+\frac{1}{\sqrt{2}}\,e^{-i\t\la}c_-(\t,\la),
\end{equation}
where  $c_{\pm}$ are the expressions between curly braces contained in \eqref{e22}.
Formula \eqref{e23} can be rewritten in a similar form
\begin{equation}\label{e25}
\wh \om_2 (\t,\la)=\frac{1}{\sqrt{2}}\,e^{i\t\la}\wh c_+(\t,\la)+\frac{1}{\sqrt{2}}\,e^{-i\t\la}\wh c_-(\t,\la),
\end{equation}
where
\begin{eqnarray}\nonumber
 &&\wh c_+(\t, \la)=
-i\Phi^*
e^{-i\t \cla^*} \S(\t)^{-1}(\cla -\la I_n)^{-1}
e^{-i\t\cla}\Phi_1
\\ &&\hspace{2cm} \times
\Big(I_{r}+
\Phi_1^* (\cla^* -\la I_n)^{-1}\Phi_2\Big),
 \label{e26}
\end{eqnarray}
\begin{eqnarray}\nonumber
 &&\wh c_-(\t, \la)=\Big(I_r+
i\Phi^*
e^{-i\t \cla^*} \S(\t)^{-1}(\cla -\la I_n)^{-1}
e^{i\t\cla}\Phi\Big)
\\ && \hspace{2cm} \times
\big(I_{r}+
\Phi^* (\cla^* -\la I_n)^{-1}\Phi_2\big).
 \label{e27}
\end{eqnarray}
In  \eqref{e26} and \eqref{e27} we used the equality $\Phi_1- \Phi=-i\Phi_2$; see the second paragraph of this section.
Comparing \eqref{e22} and \eqref{e26} yields $c_+=\wh c_+$. To prove
that $c_-=\wh c_-$ we shall need the equality
$\cla \S(\t)-\S(\t)\cla^*=i\Pi(\t)j\Pi(\t)^*$, that is, equality
(1.22) from \cite{GKSakh98} rewritten in our present notations.
Equivalently, we have
\begin{equation}\label{e28}
\S(\t)(\cla^*-\la I_n)+ie^{-i\t \cla}\Phi_1\Phi_1^*e^{i\t \cla^*}-ie^{i\t \cla}\Phi\Phi^*e^{-i\t \cla^*}=(\cla-\la I_n)\S(\t).
\end{equation}
Now, use \eqref{e22}, \eqref{e27}, and  $e^{2i\t\cla}(\Phi_1+i\Phi_2)= e^{2i\t\cla}\Phi$ to get
\begin{eqnarray}\nonumber
 &&\wh c_-(\t, \la)-c_-(\t, \la)=\Phi^*(\cla^*-\la I_n)^{-1}\Phi_2-i\Phi^*e^{-i\t \cla^*}\S(\t)^{-1}(\cla -\la I_n)^{-1}
\\  && \times \nonumber
\Big(-e^{i\t \cla}\Phi\Phi^*e^{-i\t \cla^*}+e^{-i\t \cla}\Phi_1\Phi_1^*e^{i\t \cla^*}-i\S(\t)(\cla^*-\la I_n)\Big)
\\  && \times 
e^{i\t \cla^*}(\cla^*-\la I_n)^{-1}\Phi_2.
\label{e29}
\end{eqnarray}
Finally, we  substitute \eqref{e28} into \eqref{e29}. This yields  $\wh c_-(\t, \la)=c_-(\t, \la)$.
Hence we have  $\wh c_{\pm}=c_{\pm}$, and formulas  \eqref{e24} and \eqref{e25} imply $\om_2=\wh \om_2$.
\epr

{\bf Acknowledgement.}
Daniel Alpay wishes to thank the
Earl Katz family for endowing the chair which
supported his research. 
The work of Leonid Lerer was
supported by ISF - Israel Science Foundation, Grant No~121/09, and that  of Alexander Sakhnovich  by the
Austrian Science Fund (FWF) under
Grant  No. Y330.

\vspace{2em}
\begin{flushright}
D. Alpay,\\
Department of Mathematics,\\
Ben--Gurion University of the Negev,\\
Beer-Sheva 84105,\\ Israel;\\
e-mail: {\tt dany@math.bgu.ac.il}

\vspace{1em}

M.A. Kaashoek,\\
Afdeling Wiskunde,\\
Faculteit der Exacte Wetenschappen,\\
Vrije Universiteit,\\
De Boelelaan 1081a, 1081 HV Amsterdam, 
\\
The
Netherlands;\\
e-mail: {\tt ma.kaashoek@few.vu.nl}

\vspace{1em}

L. Lerer,\\
Department of Mathematics,\\
Technion, Israel Institute of Technology,\\
Haifa 32000, Israel;\\
e-mail: {\tt llerer@techunix.technion.ac.il}

\vspace{1em}

 A.L. Sakhnovich,\\
Fakult\"at f\"ur Mathematik,\\
Universit\"at Wien,\\
Nordbergstrasse 15, 
A-1090 Wien, \\ Austria;\\
e-mail: {\tt al$\_$sakhnov@yahoo.com}

\end{flushright}

\end{document}